\documentclass[twoside]{article}

\usepackage[accepted]{aistats2025}
\usepackage{appendix}
\usepackage{amsthm}
\usepackage{amsmath}

\usepackage{amsfonts}
\usepackage{multirow}
\usepackage{multicol}

\usepackage{color}
\usepackage{algorithm}
\usepackage{algorithmic}

\usepackage{xcolor,colortbl}
\definecolor{LightCyan}{rgb}{0.88,1,1}

\usepackage{subfig}
\usepackage{tcolorbox}

\newtheorem{theorem}{Theorem}
\newtheorem{lemma}{Lemma}

\newtheorem{assumption}{Assumption}
\newtheorem{remark}{Remark}

\begin{document}

%

%

\twocolumn[

\aistatstitle{Enhanced Adaptive Gradient Algorithms for Nonconvex-PL Minimax Optimization}

\aistatsauthor{ Feihu Huang$^{1,2,*}$ \And Chunyu Xuan$^{3}$
 \And Xinrui Wang$^{1,2}$ \And Siqi Zhang$^{1,2}$ \And
 Songcan Chen$^{1,2}$}

\aistatsaddress{ 1. College of Computer Science and Technology, Nanjing University of Aeronautics and Astronautics, China; \\ 2. MIIT Key Laboratory of Pattern Analysis and Machine Intelligence, China; *huangfeihu2018@gmail.com;
\\  3. School of Automation Science and Engineering, Xi'an Jiaotong University, Xi'an, China. } ]

\begin{abstract}
  Minimax optimization recently is widely applied in many machine learning tasks such as generative adversarial networks, robust learning and reinforcement learning.
    In the paper, we study a class of nonconvex-nonconcave minimax optimization with nonsmooth regularization, where the objective function is possibly nonconvex on primal variable $x$, and it is nonconcave and satisfies the Polyak-Lojasiewicz (PL) condition on dual variable $y$. Moreover, we propose a class of enhanced momentum-based gradient descent ascent methods (i.e., MSGDA and AdaMSGDA) to solve these stochastic nonconvex-PL minimax problems. In particular, our AdaMSGDA algorithm can use various adaptive learning rates in updating the variables $x$ and $y$ without relying on any specifical types. Theoretically, we prove that our methods have the best known sample complexity of $\tilde{O}(\epsilon^{-3})$ only requiring one sample at each loop in finding an $\epsilon$-stationary solution. Some numerical experiments on PL-game and Wasserstein-GAN demonstrate the efficiency of our proposed methods.
\end{abstract}

\section{Introduction}
Minimax optimization, due to its ability of capturing the nested structures by minimizing and maximizing two subsets of variables simultaneously,
 has recently been shown great successes in many machine learning applications such as generative adversarial networks~\cite{goodfellow2014generative}, robust learning~\cite{deng2020distributionally}, AUC maximization~\cite{guo2020communication} and reinforcement learning~\cite{wai2019variance}.
In the paper, we consider the following stochastic minimax optimization problem with nonsmooth regularization
\begin{align} \label{eq:1}
 \min_{x \in \mathbb{R}^d}\max_{y\in \mathbb{R}^p} \ \mathbb{E}_{\xi \sim \mathcal{D}}\big[f(x,y;\xi)\big] + \psi(x),
\end{align}
where $f(x,y)\equiv\mathbb{E}_{\xi \sim \mathcal{D}}\big[f(x,y;\xi)\big]$ is possibly nonconvex in variable $x$, and is possibly nonconcave and satisfies Polyak-Lojasiewicz (PL) inequality~\cite{polyak1963gradient} in variable $y$. 
 Here the PL condition
is originally proposed to relax the strong convexity in minimization problem~\cite{polyak1963gradient},
which recently has been successfully analyzed in many machine learning models such as over-parameterized neural networks~\cite{frei2021proxy} and reinforcement learning~\cite{xiao2022convergence}.
In the problem~(\ref{eq:1}), the regularization $\psi(x)$ is a convex and possibly nonsmooth such as $\|x\|_1$.
Meanwhile, $\xi$ denotes a random variable drawn some
fixed but unknown distribution $\mathcal{D}$.

\begin{table*}
 \vspace*{-8pt}
  \centering
  \caption{ \textbf{Sample} or Stochastic First-Order Oracle (\textbf{SFO}) complexity comparison of the representative algorithms in finding
  an $\epsilon$-stationary solution of the \textbf{stochastic Nonconvex-PL} minimax problem \eqref{eq:1} \textbf{with or without} nonsmooth regularization $\psi(x)$, i.e., $\mathbb{E}\|\nabla F(x)\|\leq \epsilon$
  or $\mathbb{E}\|\mathcal{G}(x_t,\nabla F(x_t),\gamma)\|\leq \epsilon$. \textbf{ALR} denotes adaptive learning rate.
  Here $f(x,\cdot)$ denotes function \emph{w.r.t.} the second variable $y$ fixed $x$. \textbf{Note That} the PDAda method for Nonconvex-(PL+Concave) minimax optimization still relies on the concavity condition. Since the (Acc)SPIDER-GDA methods only focus on the finite-sum Nonconvex-PL minimax optimization ($\min_x\max_y \frac{1}{n}\sum_{i=1}^nf_i(x,y)$), we exclude them in this table.}
  \label{tab:1}
   \resizebox{\textwidth}{!}{
\begin{tabular}{c|c|c|c|c|c|c}
  \hline
   \textbf{Algorithm} & \textbf{Reference} & \textbf{Assumption } $f(x,\cdot)$ & \textbf{Complexity} & \textbf{batch size}& \textbf{ALR} & \textbf{Nonsmooth}   \\ \hline
  ZO-VRAGDA  & \cite{xu2022zeroth} & PL & $O(\epsilon^{-3})$ & $O(\epsilon^{-2})$ & $\times$ & $\times$ \\  \hline
  Smoothed-AGDA  & \cite{yang2022faster} & PL & $\tilde{O}(\epsilon^{-4})$ & $O(1)$ & $\times$ & $\times$ \\  \hline
  PDAda  & \cite{guo2021novel} & PL + Concave & $\tilde{O}(\epsilon^{-4})$ &$O(1)$ & $\checkmark$ & $\times$\\  \hline
  VRSMO & \cite{zhang2021variance} & PL  & $\tilde{O}(\epsilon^{-3})$ &$O(1)$ & $\times$ & $\times$\\  \hline
  MSGDA   & Ours & PL & {\color{red}{$\tilde{O}(\epsilon^{-3})$}} & {\color{red}{$O(1)$}} & $\times$ & {\color{red}{$\checkmark$}} \\  \hline
  AdaMSGDA   & Ours & PL & {\color{red}{$\tilde{O}(\epsilon^{-3})$}} & {\color{red}{$O(1)$}} & {\color{red}{$\checkmark$}} & {\color{red}{$\checkmark$}}  \\  \hline
\end{tabular}
 }
 \vspace*{-8pt}
\end{table*}

Due to its nested nature, the minimax problem~(\ref{eq:1}) can be rewritten as the following
minimization problem:
\begin{align} \label{eq:2}
 \min_{x\in \mathbb{R}^d} F(x) + \psi(x),
\end{align}
where $F(x) \equiv \max_{y\in \mathbb{R}^p}f(x,y)$.
Clearly, compared with the standard minimization problem, the difficulty of minimax problem~(\ref{eq:2}) is that the minimization of objective $F(x)$ depends on the maximization of objective $f(x,\cdot)$ for any $x\in \mathbb{R}^d$. To tackle this difficulty, a natural way to solve the problem~(\ref{eq:2}) by using the multi-step Gradient Descent Ascent (GDA), which uses
double-loop iterations where the inner loop uses gradient ascent to search for the approximate solution to $\max_{y\in \mathbb{R}^p}f(x,y)$ given a $x\in \mathbb{R}^d$, and the outer loop can be regarded
as running inexact gradient descent on $F(x)$ based on the approximate solution $y$.
For example, \cite{nouiehed2019solving} proposed a class of iterative first-order methods to solve nonconvex minimax problems. Subsequently, \cite{luo2020stochastic} presented a class of efficient stochastic recursive GDA to solve the stochastic Non-Convex Strongly-Concave (NC-SC) minimax problems. Meanwhile, \cite{zhang2020single} proposed an efficient single-loop smoothed GDA algorithm for nonconvex concave minimax problems.
\cite{li2022nonsmooth} presented a smoothed proximal linear descent ascent algorithm for nonsmooth composite nonconvex concave minimax optimization.
Subsequently, \cite{lu2023first} proposed a first-order augmented Lagrangian method to solve the constrained
nonconvex-concave minimax problems with nonsmooth regularization.

Another class of approaches is the alternating (two-timescale) GDA, which only uses a single-loop to update primal and dual variables $x$ and $y$ with different learning rates. For example, \cite{lin2019gradient} proposed a class of effective two-timescale (stochastic) GDA methods for nonconvex minimax optimization. Subsequently, \cite{huang2022accelerated} proposed a class of efficient momentum-based stochastic GDA methods for stochastic NC-SC minimax optimization. More recently, \cite{huang2023adagda,junchi2022nest,li2022tiada} presented some efficient adaptive stochastic GDA methods for NC-SC minimax optimization. Meanwhile,
\cite{lu2020hybrid,chen2021proximal,huang2021efficient} proposed the effective and
efficient two-timescale proximal GDA methods
to solve NC-SC minimax problems with nonsmooth regularization. In general, the two-timescale GDA methods
for minimax optimization can be more easily implemented than multi-step GDA methods,
and perform better than the multi-step GDA methods in practice.
However, the convergence analysis of the two-timescale GDA methods is more
challenging than that of the multi-step GDA methods, since the updating the primal and dual variable
$x$ and $y$ are intertwined at each step in the two-timescale GDA methods.

The above works mainly focus on the minimax optimization with a restrictive (strong) concavity condition on
the objective function $f(x,\cdot)$. In fact, many machine learning applications such as adversarial training deep neural networks (DNNs)~\cite{madry2018towards} and deep AUC maximization~\cite{guo2020communication}
do not satisfy this condition, i.e.,
the NonConvex-NonConcave (NC-NC) minimax optimization commonly appears.
Recently, some works~\cite{zhang2021variance,yang2020global,nouiehed2019solving} have begun to studying the NC-NC minimax optimization.
For example, \cite{yang2020global} proposed an alternating GDA algorithm for
minimax optimization satisfied two-sided PL condition (i.e., PL-PL), which is linearly
converges to the saddle point. \cite{nouiehed2019solving} studied the nonconvex (stochastic) minimax optimization with one-sided PL condition, and proposed a standard
multi-step GDA method to solve these Nonconvex-PL minimax problems. Subsequently, \cite{yang2022faster} developed
a class of efficient smoothed alternating GDA methods for Nonconvex-PL minimax optimization. More recently,
\cite{chen2022faster} presented a class of faster stochastic GDA (i.e.,SPIDER-GDA and
AccSPIDER-GDA) methods based on variance-reduced technique for finite-sum minimax optimization under PL condition that includes PL-PL and Nonconvex-PL minimax optimizations. Meanwhile, \cite{xu2022zeroth} proposed a class of variance-reduced zeroth-order methods for stochastic Nonconvex-PL minimax problems.

Although more recently some works have been studied the Nonconvex-PL minimax optimization, it is not well studied as the NC-SC minimax optimization, since the PL condition relaxes the strongly-convex in minimization problems (i.e., strong-concave in maximization problem)~\cite{polyak1963gradient,karimi2016linear}.
For example,
few work studies the adaptive gradient methods for Nonconvex-PL minimax optimization. Although \cite{guo2021novel} proposed an effective PDAda method for Nonconvex-(PL+Concave) minimax optimization, the PDAda method still relies on the concavity condition.

To fill this gap, thus, we propose a class of enhanced two-timescale momentum-based stochastic GDA methods (i.e., MSGDA and AdaMSGDA) to solve the stochastic Nonconvex-PL minimax problem~(\ref{eq:1}) with nonsmooth regularization. In particular, our AdaMGDA algorithm can use various types of adaptive learning rate in updating
the variables $x$ and $y$ without relying on any coordinate-wise and non-coordinate-wise
adaptive learning rates.
Specifically, our main contributions are given:
\begin{itemize}
\item[(i)] We propose a class of enhanced momentum-based gradient descent ascent algorithms (i.e., MSGDA and AdaMSGDA) to solve the Nonconvex-PL minimax problem~(\ref{eq:1}) with nonsmooth regularization. In particular, our AdaMSGDA algorithm can use various types of adaptive learning rate in updating the variables $x$ and $y$ without relying on any global and coordinate-wise adaptive learning rates.
\item[(ii)] We provide an effective convergence analysis framework for our methods. Specifically, we prove that our MSGDA and AdaMSGDA algorithms obtain the best known sample complexity of $\tilde{O}(\epsilon^{-3})$ without any large batch-sizes for finding an $\epsilon$-stationary solution of Problem~(\ref{eq:1}) (Please see Table~\ref{tab:1}), which matches a lower bound of the smooth nonconvex stochastic optimization problems~\cite{arjevani2023lower}.
\item[(iii)] Since the PL condition relaxes the strong convexity (concavity), the convergence analysis of our AdaMSGDA can be applied in the AdaGDA and VR-AdaGDA methods~\cite{huang2023adagda} for NC-SC minimax optimization. Note that \cite{huang2023adagda} only provided the convergence properties of AdaGDA and VR-AdaGDA, when they use the non-coordinate-wise adaptive learning rates to update the dual variable $y$.
\end{itemize}

\subsection*{Notations}
Given function $f(x,y)$, $f(x,\cdot)$ denotes function \emph{w.r.t.} the second variable $y$ with fixing $x$,
and $f(\cdot,y)$ denotes function \emph{w.r.t.} the first variable $x$
with fixing $y$. $I_{d}$ denotes a $d$-dimensional identity matrix. $A\succ 0$ denotes that $A \in \mathbb{R}^{d\times d}$ is positive definite,
and for any $x\in \mathbb{R}^d$ such that $\|x\|_A=\sqrt{x^TAx}$.
$\|\cdot\|$ denotes the $\ell_2$ norm for vectors and spectral norm for matrices.
$\langle x,y\rangle$ denotes the inner product of two vectors $x$ and $y$. For vectors $x$ and $y$, $x^r \ (r>0)$ denotes the element-wise
power operation, $x/y$ denotes the element-wise division and $\max(x,y)$ denotes the element-wise maximum.
$a_t=O(b_t)$ denotes that $a_t \leq C b_t$ for some constant $C>0$. The notation $\tilde{O}(\cdot)$ hides logarithmic terms.

\section{Methods}
In this section, we propose a class of enhanced two-timescale momentum-based gradient descent ascent methods (i.e., MSGDA and AdaMSGDA) to solve the stochastic Nonconvex-PL (NC-PL) minimax problem~(\ref{eq:1}) with nonsmooth regularization.

\subsection{ MSGDA Algorithm for NC-PL Minimax Optimization}
In this subsection, we consider the Nonconvex-PL stochastic minimax problem~(\ref{eq:1}), and
propose a new momentum-based stochastic gradient descent ascent (MSGDA) method.
The detailed procedure of our MSGDA method is presented in Algorithm~\ref{alg:1}.

At the lines 4 and 5 of Algorithm~\ref{alg:1}, we use the alternating proximal gradient descent ascent iterations to
update the primal variable $x$ and the dual variable $y$. Specifically, we use the proximal gradient descent to update
the primal variable $x$:
\begin{align}
&\tilde{x}_{t+1}= \mathcal{P}^{\gamma}_{\psi(\cdot)}\big(x_t - \gamma w_t\big)  \\
 & = \arg\min_{x\in \mathbb{R}^d}\Big\{
 \langle x, w_t\rangle + \frac{1}{2\gamma}(x-x_t)^TI_d(x-x_t) + \psi(x)\Big\} \nonumber \\
 & = \arg\min_{x\in \mathbb{R}^d}\Big\{
 \frac{1}{2\gamma}\|x-(x_t-\gamma w_t)\|^2 + \psi(x)\Big\}, \nonumber
\end{align}
where $\gamma$ is a constant step size and $w_t$ is a stochastic gradient estimator;
We use the gradient ascent to update
the dual variable $y$: $\tilde{y}_{t+1} = y_t + \lambda v_t$, where
$\lambda$ is a constant step size and $v_t$ is a stochastic gradient estimator.
At the lines 4 and 5 of Algorithm~\ref{alg:1}, we use the momentum iteration to further update
the primal variable $x$ and the dual variable $y$, i.e.,
\begin{align}
 y_{t+1} = y_t + \eta_t(\tilde{y}_{t+1}-y_t), \quad x_{t+1} = x_t + \eta_t(\tilde{x}_{t+1}-x_t), \nonumber
\end{align}
where $\eta_t$ is a momentum parameter shared in both the primal and dual variables.

At the lines 7 and 8 of Algorithm~\ref{alg:1}, we use momentum-based variance reduced technique (i.e., STORM~\cite{cutkosky2019momentum}) to estimate the gradients $\nabla_x f(x,y)$ and $\nabla_y f(x,y)$. Specifically, we define the stochastic gradient estimators $v_{t+1}$ and $w_{t+1}$: for all $t\geq 1$,
\begin{align}
 v_{t+1} & = \nabla_y f(x_{t+1},y_{t+1};\xi_{t+1}) \nonumber \\
 & \quad + (1-\alpha_{t+1})\big[v_t - \nabla_y f(x_t,y_t;\xi_{t+1})\big]  \\
  w_{t+1} & = \nabla_x f(x_{t+1},y_{t+1};\xi_{t+1}) \nonumber \\
 & \quad + (1-\beta_{t+1})\big[w_t - \nabla_x f(x_t,y_t;\xi_{t+1})\big],
\end{align}
where $\alpha_{t+1} \in (0,1]$ and $\beta_{t+1} \in (0,1]$ are the tuning parameters.
Note that from lines 5-8 of Algorithm~\ref{alg:1}, we simultaneously
use the momentum-based techniques to update variables and estimate the stochastic gradients.

\begin{algorithm}[t]
\caption{ Momentum-Based Gradient Descent Ascent (\textbf{MSGDA}) Algorithm}
\label{alg:1}
\begin{algorithmic}[1]
\STATE {\bfseries Input:} $T$, tuning parameters $\{\gamma, \lambda, \eta_t, \alpha_t, \beta_t\}$, initial inputs $x_1\in \mathbb{R}^d$, $y_1\in \mathbb{R}^p$; \\
\STATE {\bfseries initialize:} Draw one sample $\xi_1$,
and compute $v_1 = \nabla_y f(x_1,y_1;\xi_1)$, and $w_1 = \nabla_x f(x_1,y_1;\xi_1)$. \\
\FOR{$t=1$ \textbf{to} $T$}
\STATE $\tilde{x}_{t+1}= \mathcal{P}^{\gamma}_{\psi(\cdot)}\big(x_t - \gamma w_t\big)$, $x_{t+1} = x_t + \eta_t(\tilde{x}_{t+1}-x_t)$; \\
\STATE $\tilde{y}_{t+1} = y_t + \lambda v_t$, $y_{t+1} = y_t + \eta_t(\tilde{y}_{t+1}-y_t)$;  \\
\STATE Draw one sample $\xi_{t+1}$;
\STATE $v_{t+1} = \nabla_y f(x_{t+1},y_{t+1};\xi_{t+1}) + (1-\alpha_{t+1})\big[v_t - \nabla_y f(x_t,y_t;\xi_{t+1})\big]$; \\
\STATE $w_{t+1} = \nabla_x f(x_{t+1},y_{t+1};\xi_{t+1}) + (1-\beta_{t+1})\big[w_t - \nabla_x f(x_t,y_t;\xi_{t+1})\big] $; \\
\ENDFOR
\STATE {\bfseries Output:} Chosen uniformly random from $\{x_t,y_t\}_{t=1}^{T}$.
\end{algorithmic}
\end{algorithm}

\begin{algorithm}[t]
\caption{ Adaptive Momentum-Based Gradient Descent Ascent (\textbf{AdaMGDA}) Algorithm}
\label{alg:2}
\begin{algorithmic}[1]
\STATE {\bfseries Input:} $T$, tuning parameters $\{\gamma, \lambda, \eta_t, \alpha_t, \beta_t\}$, initial inputs $x_1\in \mathbb{R}^d$, $y_1\in \mathbb{R}^p$; \\
\STATE {\bfseries initialize:} Draw one sample $\xi_1$,
and compute $v_1 = \nabla_y f(x_1,y_1;\xi_1)$, and $w_1 = \nabla_x f(x_1,y_1;\xi_1)$. \\
\FOR{$t=1$ \textbf{to} $T$}
\STATE Generate the adaptive matrices $A_t \in \mathbb{R}^{d \times d}$ and $B_t \in \mathbb{R}^{p \times p}$;\\
\textcolor{blue}{One example of $A_t$ and $B_t$ by using update rule ($a_0 = 0$, $b_0 = 0$, $ 0 < \tau< 1$, $\rho>0$.) } \\
\textcolor{blue}{ Compute $ a_t = \tau a_{t-1} + (1 - \tau)\nabla_xf(x_t,y_t;\xi_t)^2$, $A_t = \mbox{diag}(\sqrt{a_t} + \rho)$}; \\
\textcolor{blue}{ Compute $ b_t = \tau b_{t-1} + (1 - \tau)\nabla_yf(x_t,y_t;\xi_t)^2$, $B_t = \mbox{diag}(\sqrt{b_t} + \rho)$}; \\
\STATE $\tilde{x}_{t+1}= \mathcal{P}^{\gamma}_{\psi(\cdot)}\big(x_t - \gamma A_t^{-1}w_t\big)$, $x_{t+1} = x_t + \eta_t(\tilde{x}_{t+1}-x_t)$; \\
\STATE $\tilde{y}_{t+1} = y_t + \lambda B_t^{-1}v_t$, $y_{t+1} = y_t + \eta_t(\tilde{y}_{t+1}-y_t)$;  \\
\STATE Draw one sample $\xi_{t+1}$;
\STATE $v_{t+1} = \nabla_y f(x_{t+1},y_{t+1};\xi_{t+1}) + (1-\alpha_{t+1})\big[v_t - \nabla_y f(x_t,y_t;\xi_{t+1})\big]$; \\
\STATE $w_{t+1} = \nabla_x f(x_{t+1},y_{t+1};\xi_{t+1}) + (1-\beta_{t+1})\big[w_t - \nabla_x f(x_t,y_t;\xi_{t+1})\big] $; \\
\ENDFOR
\STATE {\bfseries Output:} Chosen uniformly random from $\{x_t,y_t\}_{t=1}^{T}$.
\end{algorithmic}
\end{algorithm}

\subsection{ AdaMSGDA Algorithm for NC-PL Minimax Optimization}
In this subsection, we propose a novel adaptive momentum-based stochastic gradient descent ascent (AdaMSGDA) method to solve the Nonconvex-PL minimax problem~(\ref{eq:1}).
The detailed procedure of our AdaMSGDA method is provided in Algorithm~\ref{alg:2}.

At the lines 5 and 6 of Algorithm~\ref{alg:2}, we use the alternating adaptive gradient descent ascent iterations to
update the primal variable $x$ and the dual variable $y$. Specifically, we use adaptive stochastic
gradient descent to update
the primal variable $x$: at the $t+1$ step,
\begin{align}
 &\tilde{x}_{t+1}= \mathcal{P}^{\gamma}_{\psi(\cdot)}\big(x_t - \gamma A_t^{-1}w_t\big) \\
 & = \arg\min_{x\in \mathbb{R}^d}\Big\{
 \langle x, w_t\rangle + \frac{1}{2\gamma}(x-x_t)^TA_t(x-x_t) + \psi(x)\Big\} \nonumber \\
 & = \arg\min_{x\in \mathbb{R}^d}\Big\{
  \frac{1}{2\gamma}\|x-(x_t-\gamma A_t^{-1}w_t)\|^2_{A_t} + \psi(x)\Big\},\nonumber
\end{align}
where $\|x\|^2_A = x^TAx$, $\gamma$ is a constant step size and $w_t$ is a stochastic gradient estimator.
Here $A_t$ is an adaptive matrix representing adaptive learning rate. When $A_t$ is generated from
the one case given in Algorithm~\ref{alg:2}, as Adam algorithm~\cite{kingma2014adam}, we have
\begin{align}
 \tilde{x}_{t+1}= \mathcal{P}^{\gamma}_{\psi(\cdot)}\big(x_t - \gamma \frac{w_t}{\sqrt{a_t} + \rho}\big).
\end{align}
Meanwhile, we can also use many other forms of adaptive matrix $A_t$, e.g., we can also generate the
adaptive matrix $A_t$, as AdaBelief algorithm~\cite{zhuang2020adabelief}, defined as
\begin{align}
 & a_t = \tau a_{t-1} + (1 - \tau)(w_t-\nabla_xf(x_t,y_t;\xi_t))^2, \nonumber \\
 & A_t = \mbox{diag}(\sqrt{a_t} + \rho),
\end{align}
where $\tau\in (0,1)$.

Similarly, we use the adaptive gradient ascent to update
the dual variable $y$:
\begin{align}
\tilde{y}_{t+1} & = y_t + \lambda B_t^{-1}v_t \\
& = \arg\max_{y\in \mathbb{R}^p}\Big\{
 \langle y, v_t\rangle - \frac{1}{2\lambda}(y-y_t)^TB_t(y-y_t)\Big\}, \nonumber
\end{align}
 where $\lambda$ is a constant step size and $v_t$ is a stochastic gradient estimator.
Here $B_t$ is an adaptive matrix represented adaptive learning rate for $y$. When $B_t$ is generated from
the one case given in Algorithm~\ref{alg:2}, as Adam algorithm~\cite{kingma2014adam}, we have
\begin{align}
 \tilde{y}_{t+1}= y_t - \frac{\gamma v_t}{\sqrt{b_t+\rho}}.
\end{align}
Meanwhile, we can also use many other forms of adaptive matrix $B_t$, e.g., we can also generate the
adaptive matrix $A_t$, as AdaBelief algorithm~\cite{zhuang2020adabelief}, defined as
\begin{align}
& b_t = \tau b_{t-1} + (1 - \tau)(v_t-\nabla_yf(x_t,y_t;\xi_t))^2, \nonumber \\
& B_t = \mbox{diag}(\sqrt{b_t} + \rho).
\end{align}
where $\tau\in (0,1)$.

At the lines 5 and 6 of Algorithm~\ref{alg:2}, we use the momentum iteration to further update
the primal variable $x$ and the dual variable $y$, i.e.,
\begin{align}
 y_{t+1} = y_t + \eta_t(\tilde{y}_{t+1}-y_t), \quad x_{t+1} = x_t + \eta_t(\tilde{x}_{t+1}-x_t), \nonumber
\end{align}
where $\eta_t$ is a momentum parameter shared in both the primal and dual variables.
As Algorithm~\ref{alg:1}, at the lines 8 and 9 of Algorithm~\ref{alg:2}, we also use momentum-based variance reduced technique (i.e., STORM~\cite{cutkosky2019momentum}/ ProxHSGD~\cite{tran2022hybrid}) to estimate the gradients $\nabla_x f(x,y)$ and $\nabla_y f(x,y)$.
Note that from lines 6-9 of Algorithm~\ref{alg:2}, we also
use the momentum-based techniques to update variables and estimate the stochastic gradients simultaneously.

\section{Convergence Analysis}
In this section, we study the convergence properties of our MSGDA and AdaMSGDA algorithms
under some mild assumptions.
All related proofs are provided in the Appendix.
We first review some mild assumptions.

\begin{assumption} \label{ass:1}
Assume function $f(x,y)$ satisfies $\mu$-PL condition in variable $y\ (\mu>0)$ for
any fixed $x \in \mathbb{R}^d$, such that $\|\nabla_yf(x,y)\|^2\geq 2\mu\big(\max_{y'}f(x,y')-f(x,y)\big)$ for any $y\in \mathbb{R}^p$, where $\max_{y'}f(x,y')$ has a nonempty solution set.
\end{assumption}

\begin{assumption} \label{ass:2}
Suppose each component function $f(x,y;\xi)$ has a $L_f$-Lipschitz gradient $\nabla f(x,y;\xi)=[\nabla_x f(x,y;\xi),\nabla_y f(x,y;\xi)]$.
\end{assumption}

\begin{assumption} \label{ass:3}
Assume each component function $f(x,y;\xi)$ has an unbiased stochastic gradient $\nabla f(x,y;\xi)=[\nabla_x f(x,y;\xi),\nabla_y f(x,y;\xi)]$ with
bounded variance $\sigma^2$, i.e., $\mathbb{E}[\nabla f(x,y;\xi)] = \nabla f(x,y)$ and
\begin{align}
\mathbb{E}\|\nabla f(x,y)-\nabla f(x,y;\xi)\|^2 \leq \sigma^2. \nonumber
\end{align}
\end{assumption}

\begin{assumption} \label{ass:4}
Let $F(x)=f(x,y^*(x))=\max_y f(x,y)$.
Function $\Psi(x) = F(x)+\psi(x)$ is bounded below, \emph{i.e.,} $\Psi^* = \inf_{x\in \mathbb{R}^d}\Psi(x) > -\infty$.
\end{assumption}

Assumption~\ref{ass:1} has been commonly used in NC-PL minimax optimization~\cite{yang2022faster,chen2022faster}.
Assumptions~\ref{ass:2}-\ref{ass:3} are very commonly used in stochastic minimax optimization \cite{huang2022accelerated,xu2022zeroth}.
Assumption~\ref{ass:4} guarantees the feasibility of the minimax problem~(\ref{eq:1}).
Next, we review a useful lemma in \cite{nouiehed2019solving}.

\begin{lemma} \label{lem:AA1}
(Lemma A.5 of \cite{nouiehed2019solving})
Let $F(x)= f(x,y^*(x))=\max_y f(x,y)$ with $y^*(x) \in \
\arg\max_y f(x,y)$. Under the above Assumptions~\ref{ass:1}-\ref{ass:2},
$\nabla F(x)=\nabla_x f(x,y^*(x))$ and $F(x)$ is $L$-smooth, i.e.,
\begin{align}
\|\nabla F(x_1) - \nabla F(x_2)\| \leq L\|x_1-x_2\|, \quad \forall x_1,x_2
\end{align}
where $L=L_f(1+\frac{\kappa}{2})$ with $\kappa=\frac{L_f}{\mu}$.
\end{lemma}

\subsection{ Convergence Analysis of MSGDA Algorithm }
In this subsection, we provide the convergence properties of
our MSGDA algorithm under the above Assumptions~\ref{ass:1}-\ref{ass:4}.
We first define a useful gradient mapping
$\mathcal{G}(x,\nabla F(x),\gamma)=\frac{1}{\gamma}(x-x^+)$, where $x^+$
is generated from
\begin{align}
& x^+ = \mathcal{P}_{\psi(\cdot)}^\gamma \big(x-\gamma\nabla F(x)\big) \nonumber \\
& = \arg\min_{z\in \mathbb{R}^d}\Big\{ \langle \nabla F(x), z\rangle
+ \frac{1}{2\gamma}\|z-x\|^2 + \psi(z)\Big\}, \nonumber
\end{align}
where $\gamma>0$ and $F(x)=f(x,y^*(x))$. Then we use norm of this gradient mapping
$\|\mathcal{G}(x,\nabla F(x),\gamma)\|$ as a convergence measure, which is commonly used in
nonsmooth composite nonconvex optimization~\cite{ghadimi2016mini,j2016proximal}.

\begin{lemma} \label{lem:1}
Suppose the sequence $\{x_t,y_t\}_{t=1}^T$ be generated from Algorithm~\ref{alg:1}.
Under the Assumptions~\ref{ass:1}-\ref{ass:2}, given $0<\gamma\leq \big(\frac{\lambda\mu}{16L},\frac{\mu}{16L^2_f}\big)$ and $0<\lambda \leq \frac{1}{2L_f\eta_t}$ for all $t\geq 1$, we have
\begin{align}
& F(x_{t+1}) - f(x_{t+1},y_{t+1})
\leq (1-\frac{\eta_t\lambda\mu}{2}) \big(F(x_t) -f(x_t,y_t)\big)  \nonumber \\
& \quad + \frac{\eta_t}{8\gamma}\|\tilde{x}_{t+1}-x_t\|^2  -\frac{\eta_t}{4\lambda}\|\tilde{y}_{t+1}-y_t\|^2 \nonumber \\
& \quad + \eta_t\lambda\|\nabla_y f(x_t,y_t)-v_t\|^2,
\end{align}
where $F(x_t)=f(x_t,y^*(x_t))$ with $y^*(x_t) \in \arg\max_{y}f(x_t,y)$ for all $t\geq 1$.
\end{lemma}

Lemma~\ref{lem:1} provides the properties of the residuals $F(x_t)-f(x_t,y_t)\geq 0$ for all $t\geq 1$.
In our convergence analysis, we describe the convergence of our MSGDA algorithm by the following Lyapunov function
(i.e., potential function): for any $t\geq 1$,
\begin{align}
&\Omega_t = \mathbb{E}\big[ F(x_t) + \frac{9\gamma L^2_f}{\lambda\mu^2}\big(F(x_t) -f(x_t,y_t) \big) \nonumber +\\
& \frac{\gamma}{\eta_{t-1}}\big(\|\nabla_{x} f(x_t,y_t)-w_t\|^2 + \|\nabla_{y} f(x_t,y_t)-v_t\|^2 \big) \big]. \nonumber
\end{align}

\begin{theorem} \label{th:1}
  Under the above Assumptions~\ref{ass:1}-\ref{ass:4}, in Algorithm~\ref{alg:1}, let $\eta_t=\frac{k}{(m+t)^{1/3}}$ for all $t\geq 0$, $\alpha_{t+1}=c_1\eta_t^2$, $\beta_{t+1}=c_2\eta_t^2$, $m \geq \max\big(2,k^3, (c_1k)^3, (c_2k)^3\big)$, $k>0$, $ c_1 \geq \frac{2}{3k^3} + \frac{9L^2_f}{\mu^2}$, $c_2 \geq \frac{2}{3k^3} + \frac{9}{4}$,
 $0< \lambda \leq \min\big(\frac{3}{4\sqrt{2}\mu},\frac{m^{1/3}}{2L_fk}\big)$ and $0< \gamma \leq \big(\frac{\lambda\mu}{16L},\frac{\mu}{16L^2_f},\frac{m^{1/3}}{2Lk},\frac{1}{8L_f},\frac{\lambda\mu^2}{9L^2_f}\big)$, we have
\begin{align}
  \frac{1}{T}\sum_{t=1}^T\mathbb{E}\|\mathcal{G}(x_t,\nabla F(x_t),\gamma)\| \leq \frac{2\sqrt{3H}m^{1/6}}{T^{1/2}} + \frac{2\sqrt{3H}}{T^{1/3}}, \nonumber
\end{align}
where $\Psi(x)=F(x)+\psi(x)$ and $H = \frac{\Psi(x_1) - \Psi^*}{\gamma k} + \frac{9 L^2_f}{k\lambda\mu^2}\Delta_1 + \frac{2\sigma^2m^{1/3}}{k^2} + 2k^2(c_1^2+c_2^2)\sigma^2\ln(m+T)$ with $\Delta_1=F(x_1)-f(x_1,y_1)$.
\end{theorem}

\begin{remark}
Under the above assumptions in Theorem \ref{th:1}, we set $k=O(1)$, $c_1=O(1)$, $c_2=O(1)$,
Then we can get $m=O(1)$ and $H=\tilde{O}(1)$. Thus we have
\begin{align}
\frac{1}{T}\sum_{t=1}^T\mathbb{E}\|\mathcal{G}(x_t,\nabla F(x_t),\gamma)\| \leq \tilde{O}\Big(\frac{1}{T^{1/3}}\Big) \leq \epsilon,
\end{align}
Then we can obtain $T=\tilde{O}(\epsilon^{-3})$. Since our MSGDA algorithm requires one sample at each loop,
it has a sample complexity of $T = \tilde{O}(\epsilon^{-3})$.
Thus, our MSGDA algorithm needs $\tilde{O}(\epsilon^{-3})$
sample (or SFO) complexity for finding an $\epsilon$-stationary point.
\end{remark}

\subsection{ Convergence Analysis of AdaMSGDA Algorithm }
In this subsection, we provide the convergence properties of
our AdaMSGDA algorithm under the above Assumptions~\ref{ass:1}-\ref{ass:5}.
We first give a mild assumption~\ref{ass:5}, which is commonly used in
adaptive algorithms~\cite{huang2021super}.

\begin{assumption} \label{ass:5}
In our AdaMGDA algorithm, the adaptive matrices $A_t$ and $B_t$ for all $t\geq 1$ in updating the variables $x$ and $y$ satisfy $A_t\succeq \rho I_d \succ 0$ and $\rho_u I_{p} \succeq B_t \succeq \rho_l I_{p} \succ 0$ for any $t\geq 1$, where $\rho >0,\rho_u >0,\rho_l >0$ is an  appropriate positive number.
\end{assumption}

Assumption \ref{ass:5} ensures that the adaptive matrices $A_t$ and $B_t$ for all $t\geq 1$ are positive definite as in \cite{huang2021super,huang2023adagda}.
\textbf{Note that} \cite{huang2023adagda} only provided the convergence properties of both the AdaGDA and VR-AdaGDA when they use the global adaptive learning rates (i.e., $B_t=b_t I_p \ (b_t>0)$) to update the dual variable $y$.

Next, we define a useful gradient mapping
$\mathcal{G}(x,\nabla F(x),\gamma)=\frac{1}{\gamma}(x-x^+)$ as in~\cite{ghadimi2016mini}, where $x^+$
is generated from
\begin{align}
 & x^+ = \mathcal{P}_{\psi(\cdot)}^\gamma(x-\gamma A^{-1}\nabla F(x)) \nonumber \\
& = \arg\min_{z\in \mathbb{R}^d}\Big\{ \langle \nabla F(x), z\rangle
+ \frac{1}{2\gamma}(z-x)^TA(z-x) + \psi(z)\Big\}, \nonumber
\end{align}
where $A\succ 0$, $\gamma>0$ and $F(x)=f(x,y^*(x))$.

\begin{lemma} \label{lem:2}
Suppose the sequence $\{x_t,y_t\}_{t=1}^T$ be generated from Algorithm~\ref{alg:2}.
Under the above Assumptions, given $\gamma\leq \min\big(\frac{\lambda\mu}{16\rho_uL},\frac{\rho_l\mu}{16\rho_uL^2_f}\big)$ and $\lambda \leq \frac{1}{2\eta_t L_f\rho_u}$ for all $t\geq 1$, we have
\begin{align}
&F(x_{t+1}) - f(x_{t+1},y_{t+1})
\leq (1-\frac{\eta_t\lambda\mu}{2\rho_u}) \big(F(x_t) -f(x_t,y_t)\big) \nonumber \\
& \quad + \frac{\eta_t}{8\gamma}\|\tilde{x}_{t+1}-x_t\|^2  -\frac{\eta_t}{4\lambda\rho_u}\|\tilde{y}_{t+1}-y_t\|^2 \nonumber \\
& \quad + \frac{\eta_t\lambda}{\rho_l}\|\nabla_y f(x_t,y_t)-v_t\|^2,
\end{align}
where $F(x_t)=f(x_t,y^*(x_t))$.
\end{lemma}
Lemma~\ref{lem:2} shows the properties of the residuals $F(x_t)-f(x_t,y_t)\geq 0$ for all $t\geq 1$.
In our convergence analysis, we describe the convergence of our AdaMSGDA algorithm by the following Lyapunov function: for any $t\geq 1$,
\begin{align}
& \Phi_t = \mathbb{E}\big[ F(x_t) + \frac{9\rho_u\gamma L^2_f}{\rho\lambda\mu^2}\big(F(x_t) -f(x_t,y_t) \big) + \nonumber \\
& \frac{\gamma}{\rho\eta_{t-1}}\big(\|\nabla_{x} f(x_t,y_t)-w_t\|^2 + \|\nabla_{y} f(x_t,y_t)-v_t\|^2 \big) \big]. \nonumber
\end{align}

\begin{theorem}  \label{th:2}
Under the above Assumptions~\ref{ass:1}-\ref{ass:5}, in Algorithm~\ref{alg:2}, let $\eta_t=\frac{k}{(m+t)^{1/3}}$ for all $t\geq 0$, $\alpha_{t+1}=c_1\eta_t^2$, $\beta_{t+1}=c_2\eta_t^2$, $m \geq \max\big(2,k^3, (c_1k)^3, (c_2k)^3\big)$, $k>0$, $ c_1 \geq \frac{2}{3k^3} + \frac{9\rho_u L^2_f}{\rho_l\mu^2}$, $c_2 \geq \frac{2}{3k^3} + \frac{9}{4}$,
 $0< \lambda \leq \min\big(\frac{3}{4\sqrt{2}\mu},\frac{m^{1/3}}{2k L_f\rho_u}\big)$ and $0< \gamma \leq \min\Big(\frac{\lambda\mu}{16\rho_uL},\frac{\rho_l\mu}{16\rho_uL^2_f},\frac{m^{1/3}\rho}{2Lk},\frac{\rho}{8L_f},\frac{\rho^2\lambda\mu^2}{9\rho_uL^2_f}\Big)$, we have
\begin{align}
  \frac{1}{T}\sum_{t=1}^T\mathbb{E}\|\mathcal{G}(x_t,\nabla F(x_t),\gamma)\|  \leq \frac{2\sqrt{3G}m^{1/6}}{T^{1/2}} + \frac{2\sqrt{3G}}{T^{1/3}}, \nonumber
\end{align}
where $\Psi(x)=F(x)+\psi(x)$ and $G = \frac{\Psi(x_1) - \Psi^*}{\gamma k\rho} + \frac{9\rho_u L^2_f}{k\lambda\mu^2\rho^2}\Delta_1  + \frac{2\sigma^2m^{1/3}}{k^2\rho^2} + \frac{2k^2(c_1^2+c_2^2)\sigma^2}{\rho^2}\ln(m+T)$ and
 $\Delta_1 = F(x_1) -f(x_1,y_1)$.
\end{theorem}

\begin{figure}[ht]
\centering
 \subfloat{\includegraphics[width=0.24\textwidth]{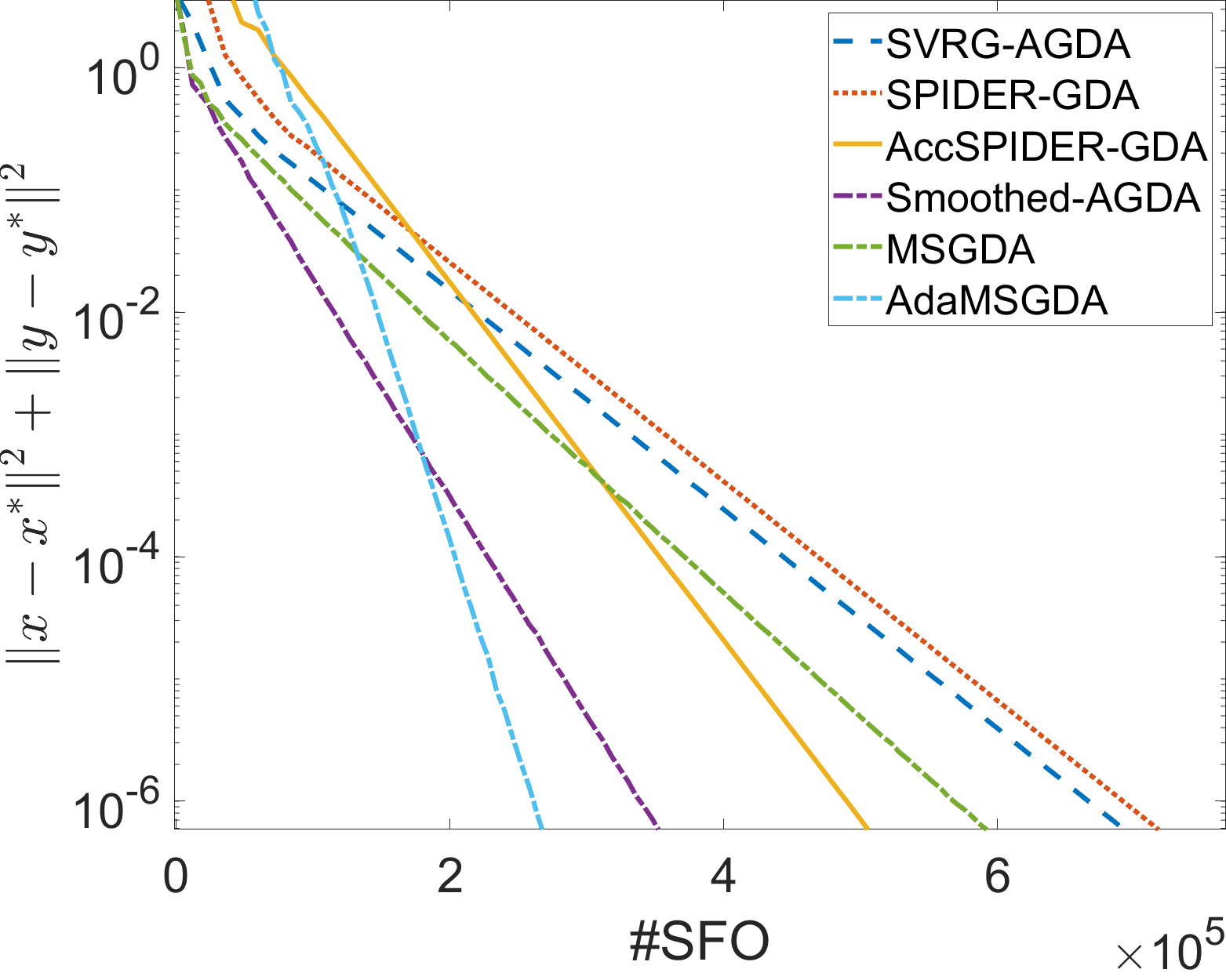}}
  \hfill
 \subfloat{\includegraphics[width=0.24\textwidth]{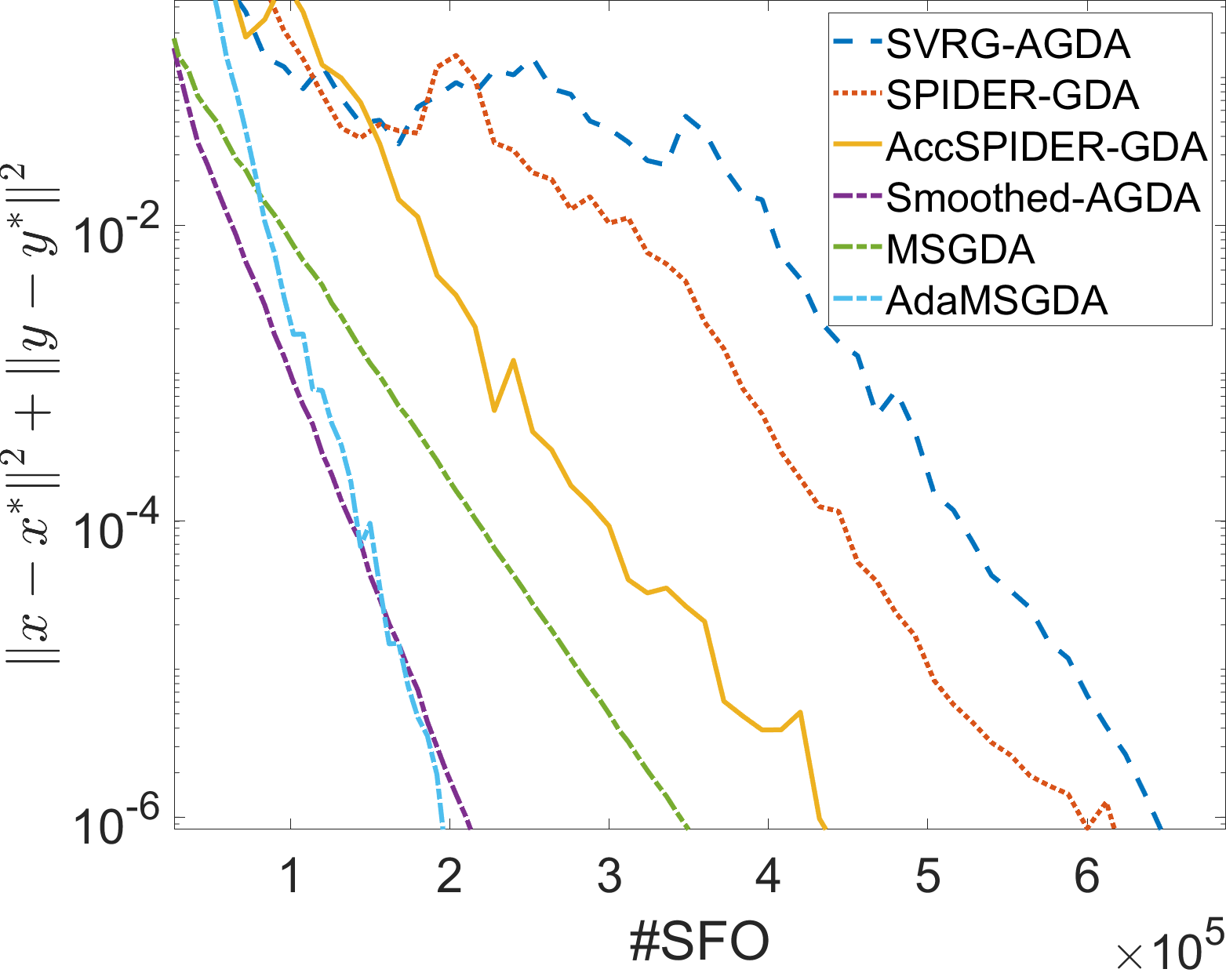}}
  \hfill
\caption{ Distance to saddle point without $L_1$ regularization:$\mu=10^{-5}$ (left), $\mu=10^{-9}$ (right).}
\label{dist}
\end{figure}

\begin{figure}[ht]
\centering
 \subfloat{\includegraphics[width=0.24\textwidth]{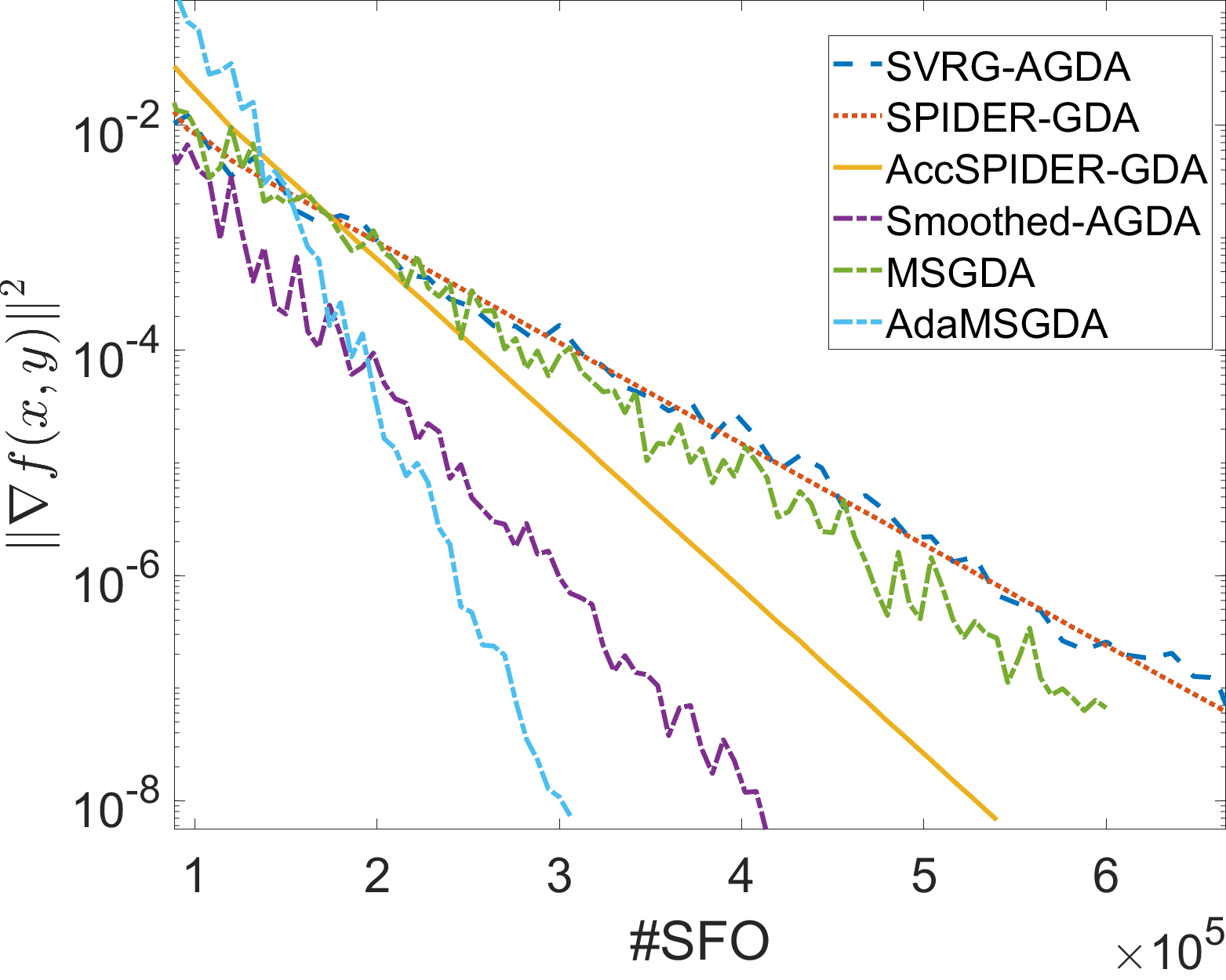}}
  \hfill
 \subfloat{\includegraphics[width=0.24\textwidth]{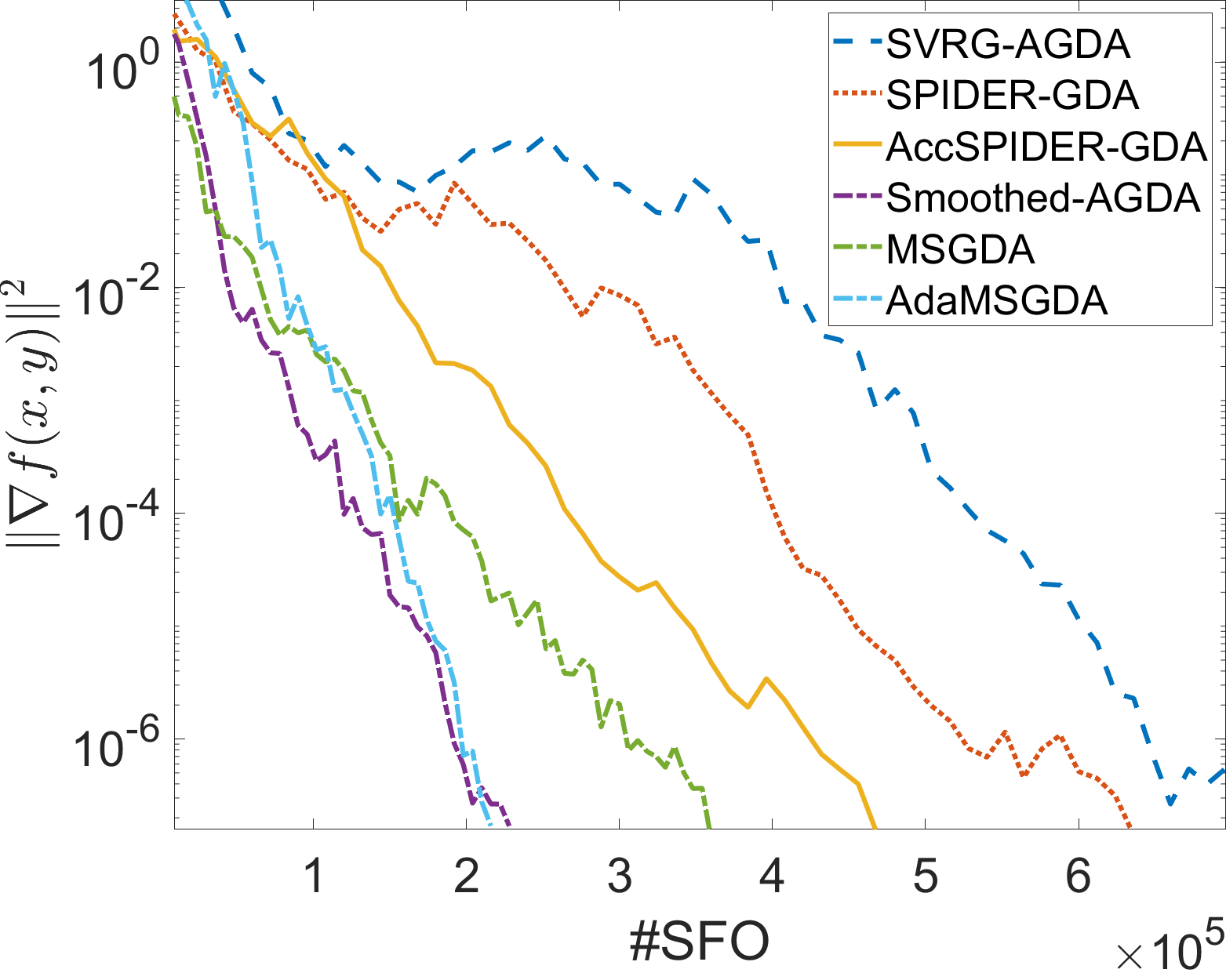}}
  \hfill
\caption{ Norm of gradient without $L_1$ regularization:$\mu=10^{-5}$ (left), $\mu=10^{-9}$ (right).}
\label{gnorm}
\end{figure}

\begin{remark}
Under the above assumptions in Theorem \ref{th:2}, we set $k=O(1)$, $\rho=\rho_l=O(1)$, $\rho_u = O(1)$ $c_1=O(1)$, $c_2=O(1)$,
we can get $m=O(1)$ and $G=\tilde{O}(1)$. Then we have
\begin{align}
\frac{1}{T}\sum_{t=1}^T\mathbb{E}\|\mathcal{G}(x_t,\nabla F(x_t),\gamma)\| \leq \tilde{O}\Big(\frac{1}{T^{1/3}}\Big) \leq \epsilon,
\end{align}
Then we can obtain $T=\tilde{O}(\epsilon^{-3})$. Since our AdaMSGDA algorithm requires one sample at each loop,
it has a sample complexity of $T = \tilde{O}(\epsilon^{-3})$.
Thus, our AdaMSGDA  algorithm needs $\tilde{O}(\epsilon^{-3})$
sample (or SFO) complexity for finding an $\epsilon$-stationary point.
\end{remark}

\section{ Numerical Experiments }
In this section, we conduct numerical experiments on Polyak-Lojasiewicz game and Wasserstein GANs to demonstrate the efficiency of our algorithms (i.e., MSGDA and AdaMSGDA). In the experiments, we compare our algorithms with the existing gradient-based minimax algorithms given in Table~\ref{tab:1}. Since the ZO-VRAGDA is a zeroth-order method, we exclude it in the comparison methods. The experiments are conducted over machine with Intel(R) Xeon(R) W-2255 CPU and Nvidia RTX2080ti(s).

\subsection{Polyak-Lojasiewicz Game}

\begin{figure}[ht]
\centering
 \subfloat{\includegraphics[width=0.24\textwidth]{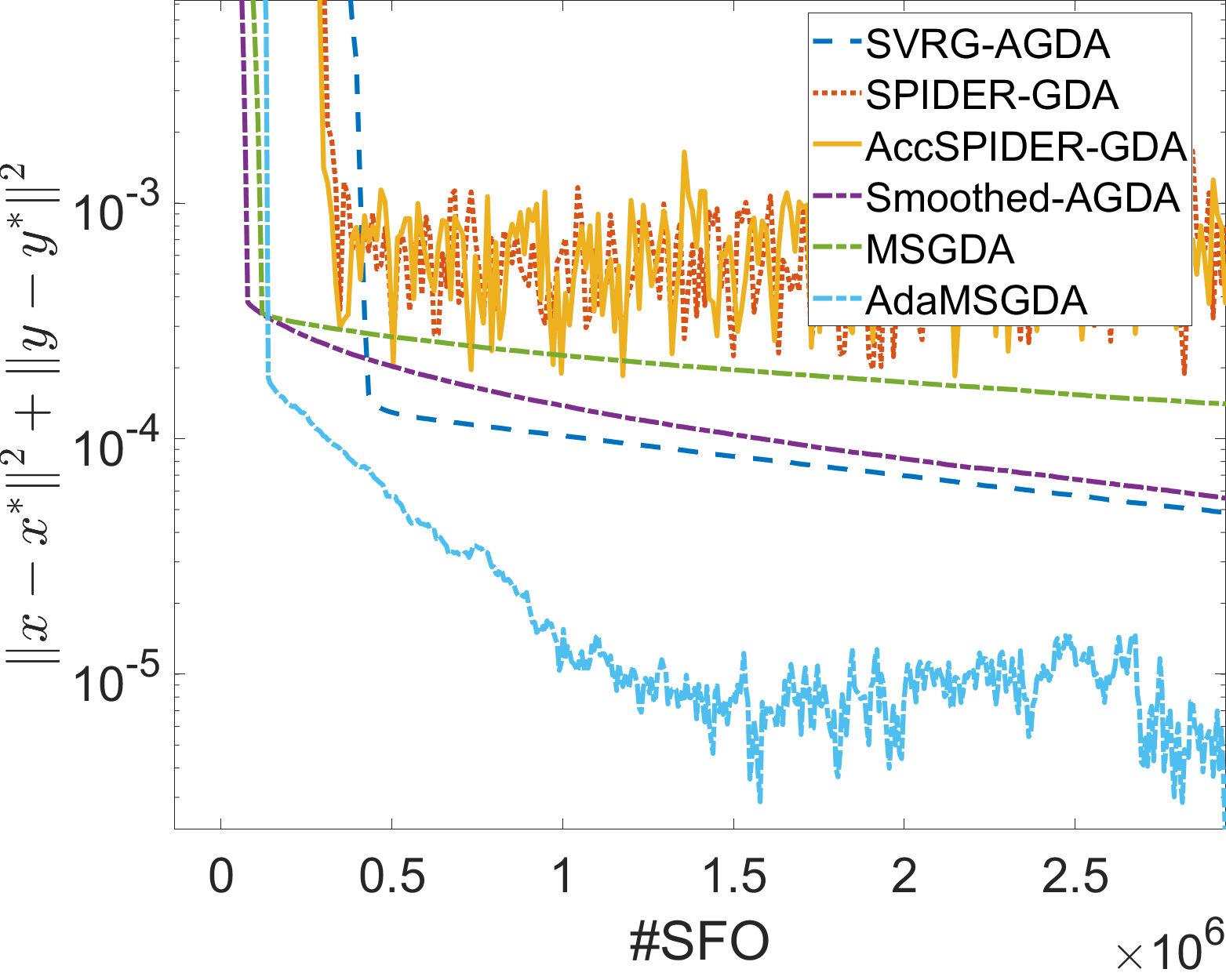}}
  \hfill
 \subfloat{\includegraphics[width=0.24\textwidth]{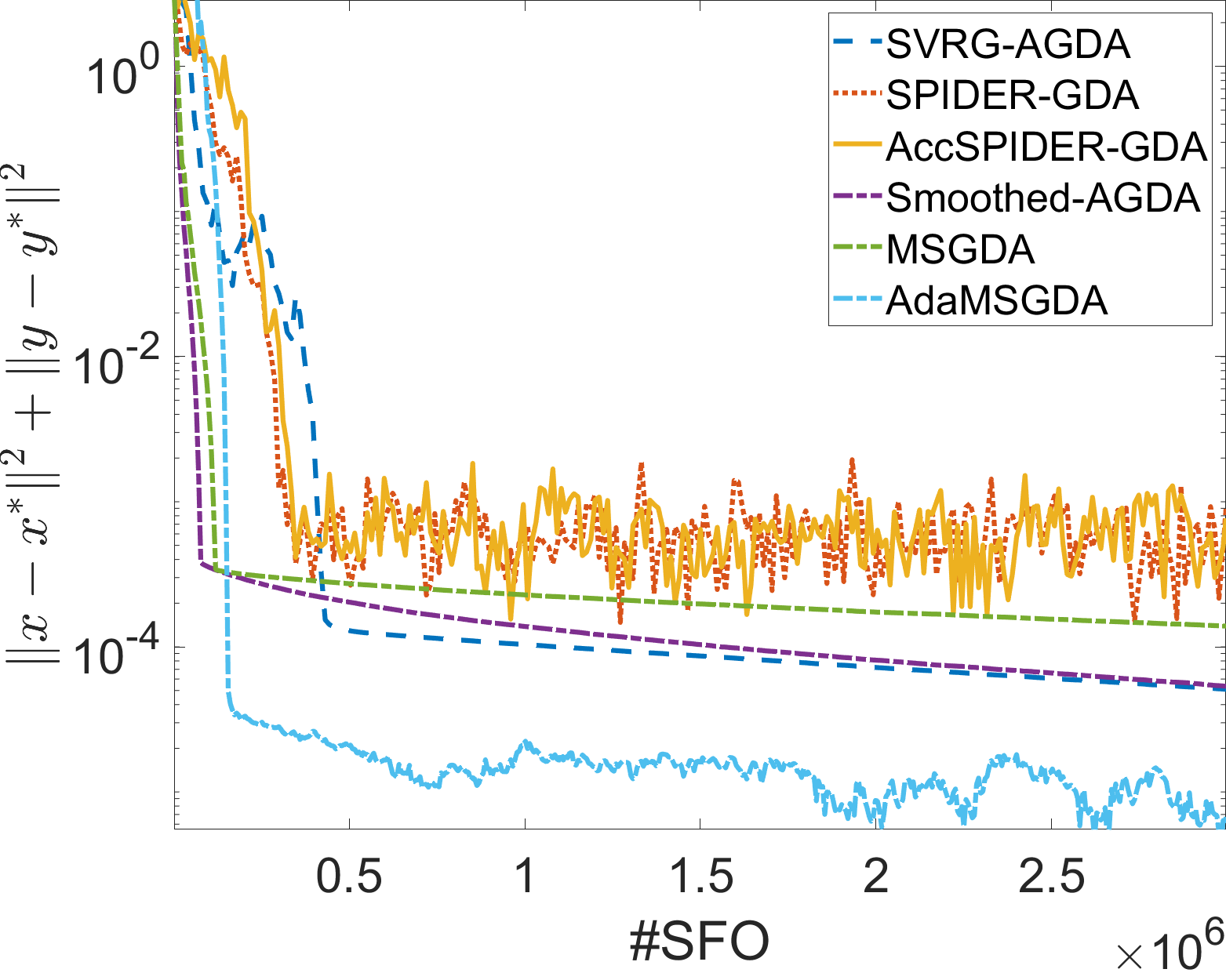}}
  \hfill
\caption{ Distance to saddle point with $L_1$  regularization:$\mu=10^{-5}$ (left), $\mu=10^{-9}$ (right).}
\label{l1dist}
\end{figure}

\begin{figure}[ht]
\centering
 \subfloat{\includegraphics[width=0.24\textwidth]{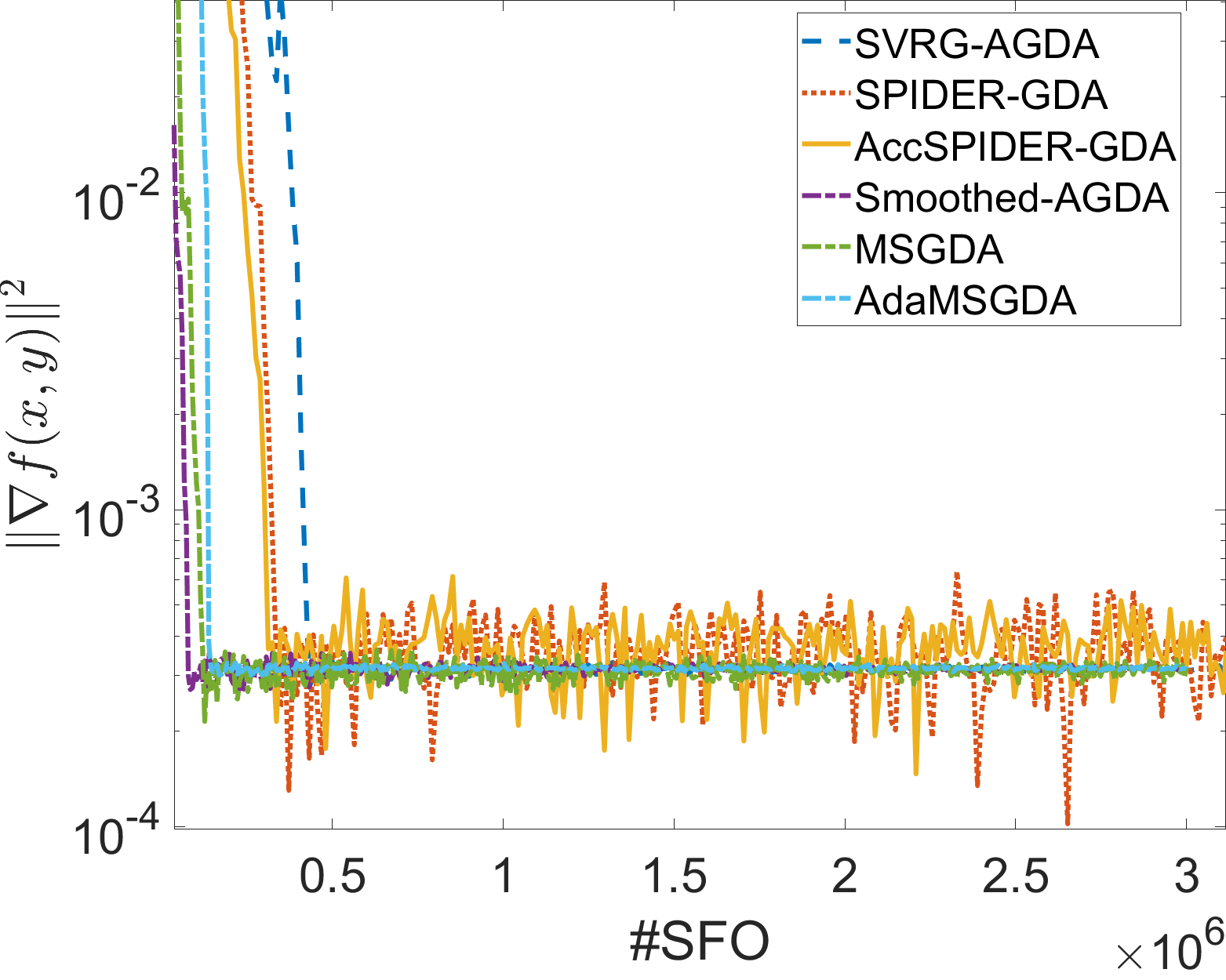}}
  \hfill
 \subfloat{\includegraphics[width=0.24\textwidth]{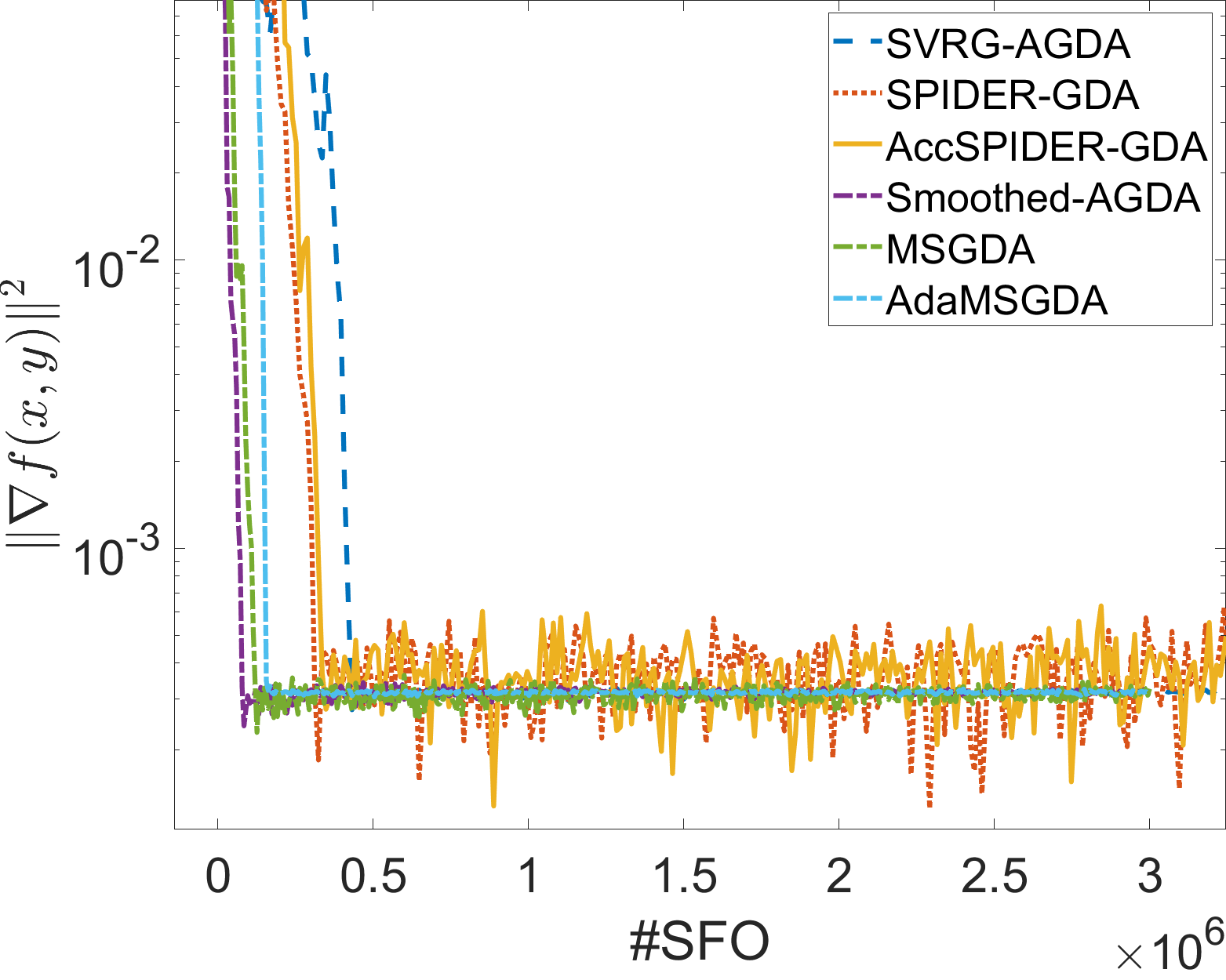}}
  \hfill
\caption{ Norm of gradient with $L_1$ regularization :$\mu=10^{-5}$ (left), $\mu=10^{-9}$ (right).}
\label{l1gnorm}
\end{figure}

In the subsection, we focus on a two-player Polyak-Lojasiewicz (PL) game, as described in \cite{chen2022faster}:
\begin{align}
    \min_{x \in \mathbb{R}^d}\max_{y\in \mathbb{R}^p} \ f(x,y)=\frac{1}{2}x^TPx-\frac{1}{2}y^TQy+x^TRy + \nu\|x\|_1, \nonumber
\end{align}
where $P=\frac{1}{n}\sum^n_{i=1}p_ip_i^T$, $Q=\frac{1}{n}\sum^n_{i=1}q_iq_i^T$ and $R=\frac{1}{n}\sum^n_{i=1}r_ir_i^T$.
Here the vectors $p_i$, $q_i$, and $r_i$ can be independently sampled from three different multivariate normal distributions: $\mathcal{N}(0,\Sigma_P)$, $\mathcal{N}(0,\Sigma_Q)$, and $\mathcal{N}(0,\Sigma_R)$, respectively. The covariance matrix $\Sigma_P$ is set in the form of $UDU^T$, where $U\in\mathbb{R}^{d\times r}$ is a column orthogonal matrix and $D\in\mathbb{R}^{r\times r}$ is a diagonal matrix with the diagonal elements distributed uniformly in the interval $[\mu, L]$, where $0<\mu < L$. 
Similarly, the matrix $\Sigma_Q$ is also set in this form. We set $\Sigma_R=0.1VV^T$, where each element of $V\in\mathbb{R}^{d\times d}$ is independently sampled from $\mathcal{N}(0,1)$. As the covariance matrices $\Sigma_P$ and $\Sigma_Q$ are rank-deficient, both matrices $P$ and $Q$ are singular. Therefore, the objective function is neither strongly convex nor strongly concave, but it satisfies the PL condition. The comparison is conducted between our algorithms and (Acc)SPIDER-GDA\cite{chen2022faster}, Smoothed-AGDA (SAGDA)\cite{yang2022faster} and the baseline SVRG-AGDA\cite{yang2020global}.
[18]. We set $n=6000$, $d=10$, $r=5$, $L=1$ for all experiments.

\begin{figure*}[ht]
    \centering
    \includegraphics[width=0.98\textwidth]{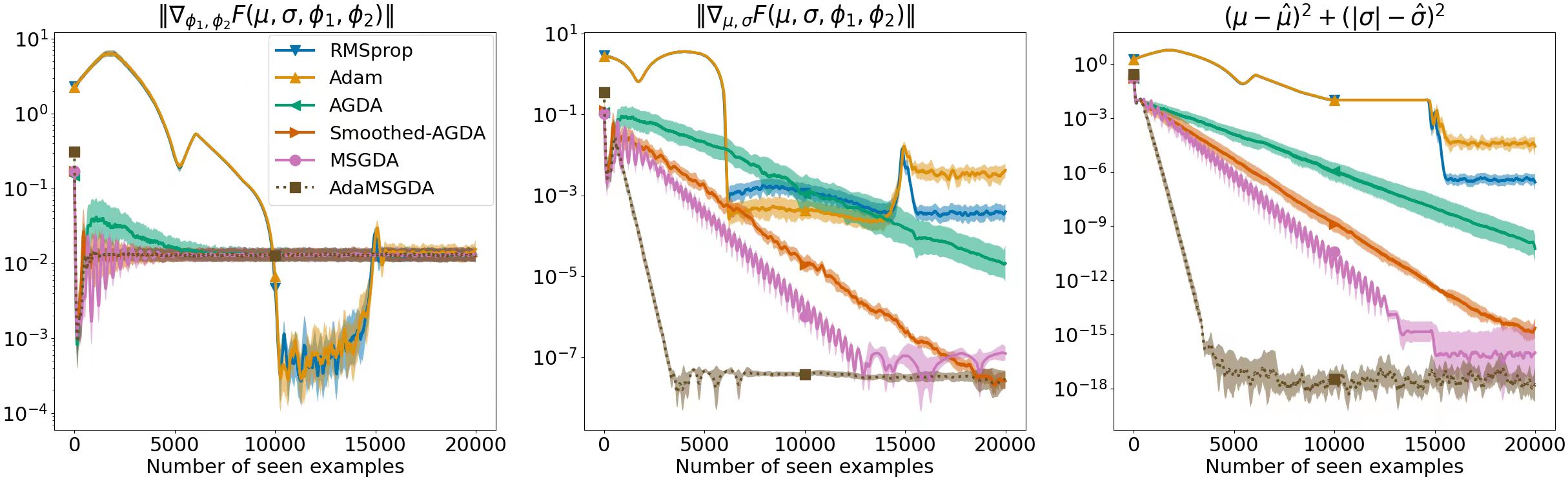}
    \caption{Training result of a Wasserstein GAN with linear generator approximating a one-dimensional Gaussian distribution}
    \label{wgan1}
\end{figure*}

\begin{figure*}[hbt]
    \centering
    \includegraphics[width=0.98\textwidth]{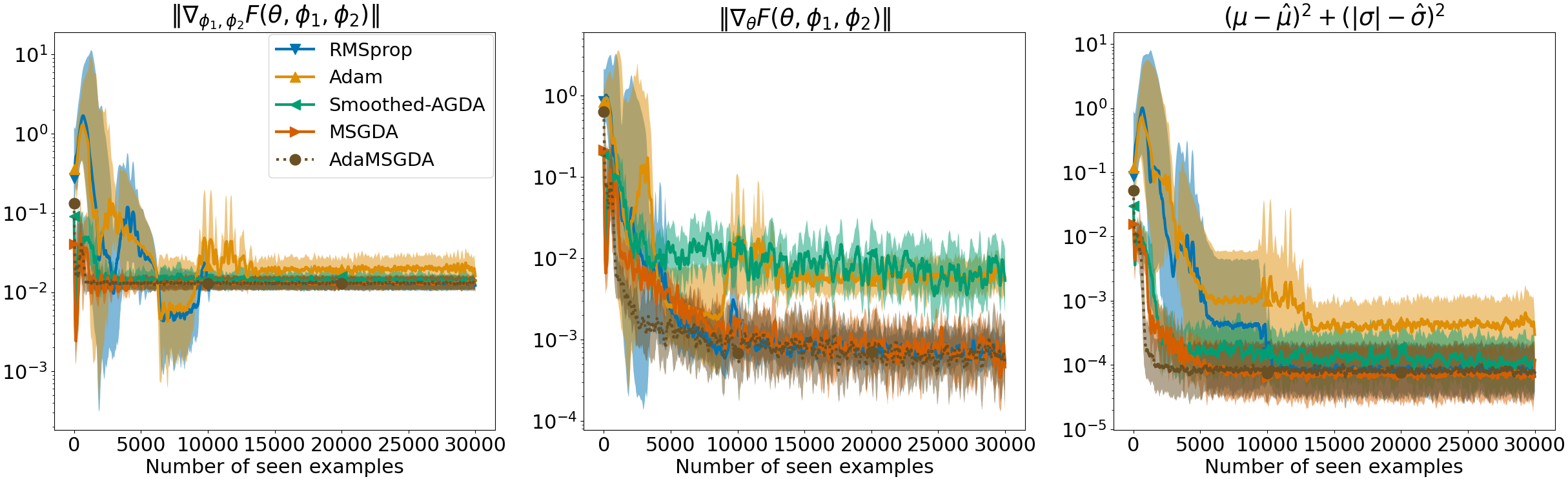}
    \caption{Training result of a Wasserstein GAN with MLP generator approximating a one-dimensional Gaussian distribution}
    \label{wgan2}
\end{figure*}

\begin{figure*}[hbt]
    \centering
    \includegraphics[width=0.98\textwidth]{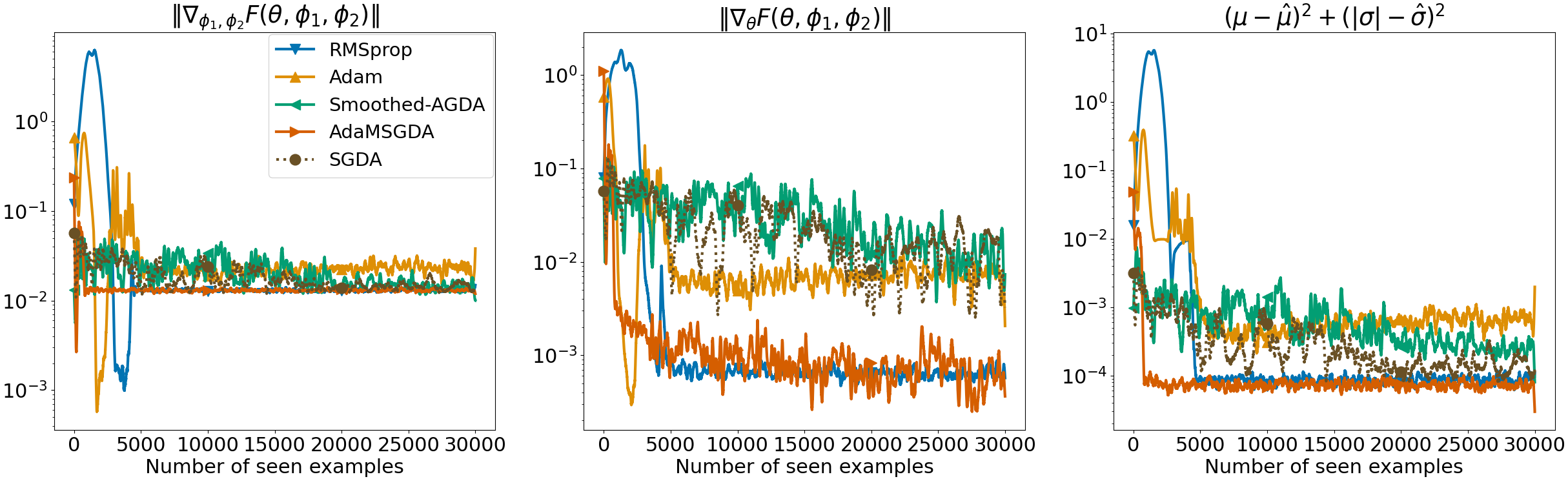}
    \caption{Training result of a Wasserstein GAN with MLP generator approximating a multi-dimensional Gaussian distribution}
    \label{wgan3}
\end{figure*}

In the experiment, we set $\eta=k/(m+t)^{1/3}$ and $\alpha=c_1\eta^2, \beta=c_2\eta^2$ for MSGDA and AdaMSGDA, where $c_1=c_2$. For AdaMGDA, we perform a grid search on $\tau$ from a range of $[0.01, 0.1, 0.2, 0.3, 0.4, 0.5, 0.7, 0.9]$ and $\rho$ from $[0.08, 0.1, 0.2, 0.3, 0.4]$ to determine the optimal parameters, which are found to be $\tau=0.01$ and $\rho=0.3$. Similarly, for Smoothed-AGDA (SAGDA), we conduct a grid search for $p$ and $\beta$, and the optimal parameters are determined to be $p=0.7$ and $\beta=0.5$. To ensure a fair comparison, we set the optimization step length for $x$, $y$, and other settings exactly the same as those used for the other methods \cite{chen2022faster}. Besides, we set $n = 6000, d = 10, r = 5, L = 1$ for all experiments.

Figures~\ref{dist}, ~\ref{gnorm}, ~\ref{l1dist} and ~\ref{l1gnorm} show the outcomes of the number of Stochastic First-Order Oracle (i.e., SFO)  calls in relation to the norm of gradient and the distance to the saddle point. As evidenced by these experimental results, our proposed methods achieve considerably faster convergence speed than existing methods under varying $\mu$ both with or without the convex nonsmooth term. As for the norm of gradient, AdaMSGDA outperforms all existing methods when $\mu=10^{-5}$ and  performs comparably to Smoothed-AGDA when $\mu=10^{-9}$.

\subsection{Wasserstein GAN}
\label{exp_gan}

In this subsection, we use a Wasserstein GAN \cite{arjovsky2017wasserstein} to approximate a one-dimensional and multi-dimensional Gaussian distribution, respectively. Specifically, we generate the data from a normal distribution $\mathcal{N}(0,0.1)$ and use a latent variable $z$ drawn from another normal distribution $\mathcal{N}(0,1)$. For a multivariate normal distribution, we generate the data from $\mathcal{N}(\mu, A)$, where $\mu = [0,0,0,0,0,0,0,0,0,0]$ and $A$ is a randomly sampled diagonal matrix with elements in $(0,1)$, ensuring it is positive-definite. Additionally, we draw a latent variable $z$ from another normal distribution $\mathcal{N}(\mu, I)$, where $I$ is the identity matrix. This Wasserstein GAN problem can be formulated as follows:
\begin{align}
& \min_{\mu,\sigma}\max_{\phi_1,\phi_2}f(\mu,\sigma,\phi_1,\phi_2) \nonumber \\
&\equiv \mathbb{E}_{x^{real},z\sim \mathcal{D}}D_{\phi}(x^{real}) -D_{\phi}(G_{\mu,\sigma}(z))-\nu||\phi||^2, \nonumber
\end{align}
where $\mathcal{D}$ is the data distribution and $\nu>0$ is the regularization parameter. Here the generator is defined as $G_{\mu,\sigma}(z)=\mu+\sigma z$, and the discriminator as $D_{\phi}(x)=\phi_1 x+\phi_2 x^2$, where $x$ denotes either real data or fake data generated by the generator. We compare our algorithms with several state-of-the-art optimization methods, including RMSprop \cite{tieleman2017divide}, Adam \cite{kingma2014adam}, AGDA \cite{yang2020global}, and Smoothed-AGDA \cite{yang2022faster}. In the experiment, we fix the batch size to 100 and repeat each algorithm three times.

In the experiment, we set $\nu=1e-3$, and fix the batch size to 100 and repeat each algorithm three times. As shown in Figure~\ref{wgan1}, the results indicate that AdaMGDA provides a significant speedup over MSGDA and achieves the best performance among all the algorithms.

Meanwhile, we consider a regularized WGAN with a neural network as a generator. For ease of comparison, we leave all the problem settings identical to the above case and only change the generator $G_{\mu,\sigma}$ to $G_\theta$, where $\theta$ are the parameters of a small neural network (one hidden layer with five neurons and ReLU activations). After tuning for each algorithm, we also
observe from Figures~\ref{wgan2},~\ref{wgan3} that our proposed MSGDA and AdaMSGDA achieve competitive performance.

In the experiment, we use the optimal parameter as \cite{yang2022faster}. Specifically, we use $\tau=5e-4,\beta_2=0.2$ for RMSprop, $\tau=5e-4,\beta_1=0.5, \beta_2=0.9$ for Adam, $\tau_1=0.1,\tau_2=0.5$ for AGDA, $\tau_1=0.1,\tau_2=0.5, \beta=0.9, P=10$ for Smoothed-AGDA and $\gamma=0.1,\lambda=0.5$ for our proposed MSGDA and AdaMSGDA. For AdaMSGDA, we use $\tau=0.01, \rho=0.3$ to form the adaptive matrix. Besides, for all the experiments we set $\eta_t=k/(m+t)^{1/3}$ and $k = 5, m = 125$.

\section{Conclusion}
In the paper, we studied a class of efficient two-timescale gradient descent ascent methods to solve nonconvex-PL minimax optimization problem with nonsmooth regularization. Specifically, we proposed an efficient momentum-based stochastic gradient descent ascent algorithm to solve the nonsmooth nonconvex-PL minimax problems. Meanwhile, we further presented an efficient
adaptive gradient-based algorithm. Moreover, we proved that our proposed methods obtain the best known and near-optimal sample complexity without relying on large batches for finding an $\epsilon$-stationary solution of nonconvex-PL minimax optimization with nonsmooth regularization.

\section*{Acknowledgements}
 This paper was partially supported by NSFC under Grant No.
62376125. It was also partially supported by the Fundamental Research Funds for the Central Universities NO.NJ2023032.

\bibliographystyle{apalike}
\bibliography{GDA-PL}

\begin{thebibliography}{}

\bibitem[Arjevani et~al., 2023]{arjevani2023lower}
Arjevani, Y., Carmon, Y., Duchi, J.~C., Foster, D.~J., Srebro, N., and
  Woodworth, B. (2023).
\newblock Lower bounds for non-convex stochastic optimization.
\newblock {\em Mathematical Programming}, 199(1-2):165--214.

\bibitem[Arjovsky et~al., 2017]{arjovsky2017wasserstein}
Arjovsky, M., Chintala, S., and Bottou, L. (2017).
\newblock Wasserstein generative adversarial networks.
\newblock In {\em International conference on machine learning}, pages
  214--223. PMLR.

\bibitem[Chen et~al., 2022]{chen2022faster}
Chen, L., Yao, B., and Luo, L. (2022).
\newblock Faster stochastic algorithms for minimax optimization under
  polyak-$\{$$\backslash$L$\}$ ojasiewicz condition.
\newblock In {\em Advances in Neural Information Processing Systems}.

\bibitem[Chen et~al., 2021]{chen2021proximal}
Chen, Z., Zhou, Y., Xu, T., and Liang, Y. (2021).
\newblock Proximal gradient descent-ascent: Variable convergence under kl
  geometry.
\newblock In {\em Proc. International Conference on Learning Representations
  (ICLR)}.

\bibitem[Cutkosky and Orabona, 2019]{cutkosky2019momentum}
Cutkosky, A. and Orabona, F. (2019).
\newblock Momentum-based variance reduction in non-convex sgd.
\newblock {\em Advances in neural information processing systems}, 32.

\bibitem[Deng et~al., 2020]{deng2020distributionally}
Deng, Y., Kamani, M.~M., and Mahdavi, M. (2020).
\newblock Distributionally robust federated averaging.
\newblock {\em Advances in neural information processing systems},
  33:15111--15122.

\bibitem[Frei and Gu, 2021]{frei2021proxy}
Frei, S. and Gu, Q. (2021).
\newblock Proxy convexity: A unified framework for the analysis of neural
  networks trained by gradient descent.
\newblock {\em Advances in Neural Information Processing Systems},
  34:7937--7949.

\bibitem[Ghadimi et~al., 2016]{ghadimi2016mini}
Ghadimi, S., Lan, G., and Zhang, H. (2016).
\newblock Mini-batch stochastic approximation methods for nonconvex stochastic
  composite optimization.
\newblock {\em Mathematical Programming}, 155(1-2):267--305.

\bibitem[Goodfellow et~al., 2014]{goodfellow2014generative}
Goodfellow, I., Pouget-Abadie, J., Mirza, M., Xu, B., Warde-Farley, D., Ozair,
  S., Courville, A., and Bengio, Y. (2014).
\newblock Generative adversarial nets.
\newblock In {\em Advances in neural information processing systems}, pages
  2672--2680.

\bibitem[Guo et~al., 2020]{guo2020communication}
Guo, Z., Liu, M., Yuan, Z., Shen, L., Liu, W., and Yang, T. (2020).
\newblock Communication-efficient distributed stochastic auc maximization with
  deep neural networks.
\newblock In {\em International conference on machine learning}, pages
  3864--3874. PMLR.

\bibitem[Guo et~al., 2021]{guo2021novel}
Guo, Z., Xu, Y., Yin, W., Jin, R., and Yang, T. (2021).
\newblock A novel convergence analysis for algorithms of the adam family and
  beyond.
\newblock {\em arXiv preprint arXiv:2104.14840}.

\bibitem[Huang et~al., 2022]{huang2022accelerated}
Huang, F., Gao, S., Pei, J., and Huang, H. (2022).
\newblock Accelerated zeroth-order and first-order momentum methods from mini
  to minimax optimization.
\newblock {\em Journal of Machine Learning Research}, 23(36):1--70.

\bibitem[Huang et~al., 2021a]{huang2021super}
Huang, F., Li, J., and Huang, H. (2021a).
\newblock Super-adam: faster and universal framework of adaptive gradients.
\newblock {\em Advances in Neural Information Processing Systems},
  34:9074--9085.

\bibitem[Huang et~al., 2023]{huang2023adagda}
Huang, F., Wu, X., and Hu, Z. (2023).
\newblock Adagda: Faster adaptive gradient descent ascent methods for minimax
  optimization.
\newblock In {\em International Conference on Artificial Intelligence and
  Statistics}, pages 2365--2389. PMLR.

\bibitem[Huang et~al., 2021b]{huang2021efficient}
Huang, F., Wu, X., and Huang, H. (2021b).
\newblock Efficient mirror descent ascent methods for nonsmooth minimax
  problems.
\newblock {\em Advances in Neural Information Processing Systems},
  34:10431--10443.

\bibitem[J~Reddi et~al., 2016]{j2016proximal}
J~Reddi, S., Sra, S., Poczos, B., and Smola, A.~J. (2016).
\newblock Proximal stochastic methods for nonsmooth nonconvex finite-sum
  optimization.
\newblock {\em Advances in neural information processing systems}, 29.

\bibitem[Junchi et~al., 2022]{junchi2022nest}
Junchi, Y., Li, X., and He, N. (2022).
\newblock Nest your adaptive algorithm for parameter-agnostic nonconvex minimax
  optimization.
\newblock In {\em Advances in Neural Information Processing Systems}.

\bibitem[Karimi et~al., 2016]{karimi2016linear}
Karimi, H., Nutini, J., and Schmidt, M. (2016).
\newblock Linear convergence of gradient and proximal-gradient methods under
  the polyak-{\l}ojasiewicz condition.
\newblock In {\em Joint European conference on machine learning and knowledge
  discovery in databases}, pages 795--811. Springer.

\bibitem[Kingma and Ba, 2014]{kingma2014adam}
Kingma, D.~P. and Ba, J. (2014).
\newblock Adam: A method for stochastic optimization.
\newblock {\em arXiv preprint arXiv:1412.6980}.

\bibitem[Li et~al., 2022a]{li2022nonsmooth}
Li, J., Zhu, L., and So, A. M.-C. (2022a).
\newblock Nonsmooth composite nonconvex-concave minimax optimization.
\newblock {\em arXiv preprint arXiv:2209.10825}.

\bibitem[Li et~al., 2022b]{li2022tiada}
Li, X., Yang, J., and He, N. (2022b).
\newblock Tiada: A time-scale adaptive algorithm for nonconvex minimax
  optimization.
\newblock {\em arXiv preprint arXiv:2210.17478}.

\bibitem[Lin et~al., 2020]{lin2019gradient}
Lin, T., Jin, C., and Jordan, M. (2020).
\newblock On gradient descent ascent for nonconvex-concave minimax problems.
\newblock In {\em International Conference on Machine Learning}, pages
  6083--6093. PMLR.

\bibitem[Lu et~al., 2020]{lu2020hybrid}
Lu, S., Tsaknakis, I., Hong, M., and Chen, Y. (2020).
\newblock Hybrid block successive approximation for one-sided non-convex
  min-max problems: algorithms and applications.
\newblock {\em IEEE Transactions on Signal Processing}, 68:3676--3691.

\bibitem[Lu and Mei, 2023]{lu2023first}
Lu, Z. and Mei, S. (2023).
\newblock A first-order augmented lagrangian method for constrained minimax
  optimization.
\newblock {\em arXiv preprint arXiv:2301.02060}.

\bibitem[Luo et~al., 2020]{luo2020stochastic}
Luo, L., Ye, H., Huang, Z., and Zhang, T. (2020).
\newblock Stochastic recursive gradient descent ascent for stochastic
  nonconvex-strongly-concave minimax problems.
\newblock {\em Advances in Neural Information Processing Systems}, 33.

\bibitem[Madry et~al., 2018]{madry2018towards}
Madry, A., Makelov, A., Schmidt, L., Tsipras, D., and Vladu, A. (2018).
\newblock Towards deep learning models resistant to adversarial attacks.

\bibitem[Nouiehed et~al., 2019]{nouiehed2019solving}
Nouiehed, M., Sanjabi, M., Huang, T., Lee, J.~D., and Razaviyayn, M. (2019).
\newblock Solving a class of non-convex min-max games using iterative first
  order methods.
\newblock {\em Advances in Neural Information Processing Systems}, 32.

\bibitem[Polyak, 1963]{polyak1963gradient}
Polyak, B. (1963).
\newblock Gradient methods for the minimisation of functionals.
\newblock {\em USSR Computational Mathematics and Mathematical Physics},
  3(4):864--878.

\bibitem[Tieleman and Hinton, 2017]{tieleman2017divide}
Tieleman, T. and Hinton, G. (2017).
\newblock Divide the gradient by a running average of its recent magnitude.
  coursera: Neural networks for machine learning.
\newblock {\em Technical report}.

\bibitem[Tran-Dinh et~al., 2022]{tran2022hybrid}
Tran-Dinh, Q., Pham, N.~H., Phan, D.~T., and Nguyen, L.~M. (2022).
\newblock A hybrid stochastic optimization framework for composite nonconvex
  optimization.
\newblock {\em Mathematical Programming}, 191(2):1005--1071.

\bibitem[Wai et~al., 2019]{wai2019variance}
Wai, H.-T., Hong, M., Yang, Z., Wang, Z., and Tang, K. (2019).
\newblock Variance reduced policy evaluation with smooth function
  approximation.
\newblock {\em Advances in Neural Information Processing Systems},
  32:5784--5795.

\bibitem[Xiao, 2022]{xiao2022convergence}
Xiao, L. (2022).
\newblock On the convergence rates of policy gradient methods.
\newblock {\em Journal of Machine Learning Research}, 23(282):1--36.

\bibitem[Xu et~al., 2022]{xu2022zeroth}
Xu, Z., Wang, Z.-Q., Wang, J.-L., and Dai, Y.-H. (2022).
\newblock Zeroth-order alternating gradient descent ascent algorithms for a
  class of nonconvex-nonconcave minimax problems.
\newblock {\em arXiv preprint arXiv:2211.13668}.

\bibitem[Yang et~al., 2020]{yang2020global}
Yang, J., Kiyavash, N., and He, N. (2020).
\newblock Global convergence and variance reduction for a class of
  nonconvex-nonconcave minimax problems.
\newblock {\em Advances in Neural Information Processing Systems}, 33.

\bibitem[Yang et~al., 2022]{yang2022faster}
Yang, J., Orvieto, A., Lucchi, A., and He, N. (2022).
\newblock Faster single-loop algorithms for minimax optimization without strong
  concavity.
\newblock In {\em International Conference on Artificial Intelligence and
  Statistics}, pages 5485--5517. PMLR.

\bibitem[Zhang et~al., 2020]{zhang2020single}
Zhang, J., Xiao, P., Sun, R., and Luo, Z. (2020).
\newblock A single-loop smoothed gradient descent-ascent algorithm for
  nonconvex-concave min-max problems.
\newblock {\em Advances in neural information processing systems},
  33:7377--7389.

\bibitem[Zhang, 2021]{zhang2021variance}
Zhang, L. (2021).
\newblock Variance reduction for non-convex stochastic optimization: General
  analysis and new applications.
\newblock Master's thesis, ETH Zurich.

\bibitem[Zhuang et~al., 2020]{zhuang2020adabelief}
Zhuang, J., Tang, T., Ding, Y., Tatikonda, S.~C., Dvornek, N., Papademetris,
  X., and Duncan, J. (2020).
\newblock Adabelief optimizer: Adapting stepsizes by the belief in observed
  gradients.
\newblock {\em Advances in neural information processing systems},
  33:18795--18806.

\end{thebibliography}

\section*{Checklist}

%

 \begin{enumerate}

 \item For all models and algorithms presented, check if you include:
 \begin{enumerate}
   \item A clear description of the mathematical setting, assumptions, algorithm, and/or model. [Yes]
   \item An analysis of the properties and complexity (time, space, sample size) of any algorithm. [Yes]
   \item (Optional) Anonymized source code, with specification of all dependencies, including external libraries. [Not Applicable]
 \end{enumerate}

 \item For any theoretical claim, check if you include:
 \begin{enumerate}
   \item Statements of the full set of assumptions of all theoretical results. [Yes]
   \item Complete proofs of all theoretical results. [Yes]
   \item Clear explanations of any assumptions. [Yes]
 \end{enumerate}

 \item For all figures and tables that present empirical results, check if you include:
 \begin{enumerate}
   \item The code, data, and instructions needed to reproduce the main experimental results (either in the supplemental material or as a URL). [Yes]
   \item All the training details (e.g., data splits, hyperparameters, how they were chosen). [Yes]
         \item A clear definition of the specific measure or statistics and error bars (e.g., with respect to the random seed after running experiments multiple times). [Yes]
         \item A description of the computing infrastructure used. (e.g., type of GPUs, internal cluster, or cloud provider). [Yes]
 \end{enumerate}

 \item If you are using existing assets (e.g., code, data, models) or curating/releasing new assets, check if you include:
 \begin{enumerate}
   \item Citations of the creator If your work uses existing assets. [Yes]
   \item The license information of the assets, if applicable. [Not Applicable]
   \item New assets either in the supplemental material or as a URL, if applicable. [Yes]
   \item Information about consent from data providers/curators. [Not Applicable]
   \item Discussion of sensible content if applicable, e.g., personally identifiable information or offensive content. [Not Applicable]
 \end{enumerate}

 \item If you used crowdsourcing or conducted research with human subjects, check if you include:
 \begin{enumerate}
   \item The full text of instructions given to participants and screenshots. [Not Applicable]
   \item Descriptions of potential participant risks, with links to Institutional Review Board (IRB) approvals if applicable. [Not Applicable]
   \item The estimated hourly wage paid to participants and the total amount spent on participant compensation. [Not Applicable]
 \end{enumerate}

 \end{enumerate}

\begin{onecolumn}

\appendix

\section{Appendix}
In this section, we provide the detailed convergence analysis of our algorithms.
We first review some useful lemmas.

\begin{lemma} \label{lem:A1}
(Lemma A.5 of \cite{nouiehed2019solving})
Let $F(x)= f(x,y^*(x))=\max_y f(x,y)$ with $y^*(x) \in \
\arg\max_y f(x,y)$. Under the above Assumptions~\ref{ass:1}-\ref{ass:2},
$\nabla F(x)=\nabla_x f(x,y^*(x))$ and $F(x)$ is $L$-smooth, i.e.,
\begin{align}
\|\nabla F(x_1) - \nabla F(x_2)\| \leq L\|x_1-x_2\|, \quad \forall x_1,x_2
\end{align}
where $L=L_f(1+\frac{\kappa}{2})$ with $\kappa=\frac{L_f}{\mu}$.
\end{lemma}

\begin{lemma} \label{lem:A2}
(\cite{karimi2016linear})
 Function $f(x): \mathbb{R}^d\rightarrow \mathbb{R}$ is $L$-smooth and satisfies PL condition with constant $\mu$, then it also
satisfies error bound (EB) condition with $\mu$, i.e., for all $x \in \mathbb{R}^d$
\begin{align}
 \|\nabla f(x)\| \geq \mu\|x^*-x\|,
\end{align}
where $x^* \in \arg\min_{x} f(x)$. It also satisfies quadratic growth (QG) condition with $\mu$, i.e.,
\begin{align}
 f(x)-\min_x f(x) \geq \frac{\mu}{2}\|x^*-x\|^2.
\end{align}
\end{lemma}
From the above lemma~\ref{lem:A2}, when consider the problem $\max_x f(x)$ that is equivalent to the problem $-\min_x -f(x)$, we have
\begin{align}
 & \|\nabla f(x)\| \geq \mu\|x^*-x\|, \\
 & \max_x f(x) - f(x) \geq \frac{\mu}{2}\|x^*-x\|^2.
\end{align}

\begin{lemma} \label{lem:B1}
Under the above Assumptions~\ref{ass:1}-\ref{ass:2}, and assume the gradient estimators
$\big\{v_t,w_t\big\}_{t=1}^T$ be generated from Algorithm \ref{alg:1},
we have
 \begin{align}
 \mathbb{E}\|v_{t+1} - \nabla_y f(x_{t+1},y_{t+1})\|^2
 & \leq (1-\alpha_{t+1})\mathbb{E} \|v_t - \nabla_y f(x_t,y_t)\|^2 + 2\alpha_{t+1}^2\sigma^2 \nonumber \\
 & \quad + 4L_f^2\eta_t^2\mathbb{E}\big(\|\tilde{x}_{t+1} - x_t\|^2 + \|\tilde{y}_{t+1} - y_t\|^2 \big),
 \end{align}
 \begin{align}
\mathbb{E}\|w_{t+1} - \nabla_x f(x_{t+1},y_{t+1}) \|^2 & \leq (1-\beta_{t+1}) \mathbb{E}\|w_t - \nabla_x f(x_t,y_t)\|^2 + 2\beta^2_{t+1}\sigma^2 \nonumber \\
& \quad + 4L^2_f\eta^2_t\mathbb{E}\big(\|\tilde{x}_{t+1}-x_t\|^2 + \|\tilde{y}_{t+1}-y_t\|^2 \big).
\end{align}
\end{lemma}

\begin{proof}
Since $w_{t+1} = \nabla_x f(x_{t+1},y_{t+1};\xi_{t+1}) + (1-\beta_{t+1})\big(w_t
 - \nabla_x f(x_t,y_t;\xi_{t+1})\big)$,
we have
\begin{align}
 &\mathbb{E}\|w_{t+1} - \nabla_x f(x_{t+1},y_{t+1})\|^2  \\
 & = \mathbb{E}\big\|\nabla_x f(x_{t+1},y_{t+1};\xi_{t+1}) + (1-\beta_{t+1})\big(w_t
 - \nabla_x f(x_t,y_t;\xi_{t+1})\big) - \nabla_x f(x_{t+1},y_{t+1})\big\|^2 \nonumber \\
 & = \mathbb{E}\big\|\nabla_x f(x_{t+1},y_{t+1};\xi_{t+1}) - \nabla_x f(x_{t+1},y_{t+1}) - (1-\beta_{t+1})\big(
 \nabla_x f(x_t,y_t;\xi_{t+1}) - \nabla_x f(x_t,y_t) \big) \nonumber \\
 & \quad + (1-\beta_{t+1})\big( w_t - \nabla_x f(x_t,y_t) \big) \big\|^2 \nonumber \\
 & \mathop{=}^{(i)} \mathbb{E}\big\|\nabla_x f(x_{t+1},y_{t+1};\xi_{t+1}) - \nabla_xf(x_{t+1},y_{t+1})  - (1-\beta_{t+1})\big(
 \nabla_x f(x_t,y_t;\xi_{t+1}) - \nabla_x f(x_t,y_t) \big) \big\|^2 \nonumber \\
 & \quad + (1-\beta_{t+1})^2\mathbb{E}\|w_t - \nabla_x f(x_t,y_t)\|^2 \nonumber \\
 & \leq 2(1-\beta_{t+1})^2\mathbb{E}\big\|\nabla_x f(x_{t+1},y_{t+1};\xi_{t+1}) - \nabla_x f(x_t,y_t;\xi_{t+1}) -
 \big( \nabla_x f(x_{t+1},y_{t+1})- \nabla_x f(x_t,y_t)\big) \big\|^2 \nonumber \\
 & \quad + 2\beta^2_{t+1}\mathbb{E}\big\|\nabla_x f(x_{t+1},y_{t+1};\xi_{t+1}) -
 \nabla_x f(x_{t+1},y_{t+1}) \big\|^2 + (1-\beta_{t+1})^2\mathbb{E}\|w_t - \nabla_x f(x_t,y_t)\|^2 \nonumber \\
 & \mathop{\leq}^{(ii)} 2(1-\beta_{t+1})^2\mathbb{E}\big\|\nabla_x f(x_{t+1},y_{t+1};\xi_{t+1}) - \nabla_x f(x_t,y_t;\xi_{t+1}) \big\|^2 + 2\beta^2_{t+1}\sigma^2 \nonumber \\
 & \quad  + (1-\beta_{t+1})^2\mathbb{E}\|w_t - \nabla_x f(x_t,y_t)\|^2 \nonumber \\
 & \leq (1-\beta_{t+1})^2 \mathbb{E}\|w_t - \nabla_x f(x_t,y_t)\|^2 + 2\beta^2_{t+1}\sigma^2 + 4(1-\beta_{t+1})^2L^2_f \mathbb{E}\big(
  \|x_{t+1}-x_t\|^2 + \|y_{t+1}-y_t\|^2 \big) \nonumber \\
 & \leq (1-\beta_{t+1}) \mathbb{E}\|w_t - \nabla_x f(x_t,y_t)\|^2 + 2\beta^2_{t+1}\sigma^2 + 4L^2_f\eta^2_t\mathbb{E}\big(
  \|\tilde{x}_{t+1}-x_t\|^2 + \|\tilde{y}_{t+1}-y_t\|^2 \big) \nonumber,
\end{align}
where the equality (i) holds by
\begin{align}
 \mathbb{E}_{\xi_{t+1}} \big[\nabla_x f(x_{t+1},y_{t+1};\xi_{t+1})-\nabla_x f(x_{t+1},y_{t+1})\big] = 0, \
 \mathbb{E}_{\xi_{t+1}} \big[\nabla_x f(x_t,y_t;\xi_{t+1})) - \nabla_x f(x_t,y_t) \big] = 0, \nonumber
\end{align}
and the inequality (ii) holds by the inequality $\mathbb{E}\|\zeta-\mathbb{E}[\zeta]\|^2 \leq \mathbb{E}\|\zeta\|^2$ and Assumption \ref{ass:3};
the second last inequality is due to Assumption \ref{ass:2}; the last inequality holds by $0<\beta_{t+1} \leq 1$ and $x_{t+1}=x_t+\eta_t(\tilde{x}_{t+1}-x_t)$,
$y_{t+1}=y_t+\eta_t(\tilde{y}_{t+1}-y_t)$.

Similarly, we can obtain
 \begin{align}
 \mathbb{E}\|v_{t+1} - \nabla_y f(x_{t+1},y_{t+1})\|^2
 & \leq (1-\alpha_{t+1})\mathbb{E} \|v_t - \nabla_y f(x_t,y_t)\|^2 + 2\alpha_{t+1}^2\sigma^2 \nonumber \\
 & \quad + 4L_f^2\eta_t^2\mathbb{E}\big(\|\tilde{x}_{t+1} - x_t\|^2 + \|\tilde{y}_{t+1} - y_t\|^2 \big).
 \end{align}

\end{proof}

\subsection{ Convergence Analysis of MSGDA Algorithm }
In this subsection, we detail the convergence analysis of MSGDA algorithm.

\begin{lemma} \label{lem:C1}
(Restatement of Lemma 2)
Suppose the sequence $\{x_t,y_t\}_{t=1}^T$ be generated from Algorithm~\ref{alg:1}.
Under the Assumptions~\ref{ass:1}-\ref{ass:2}, given $\gamma\leq \big(\frac{\lambda\mu}{16L},\frac{\mu}{16L^2_f}\big)$ and $0<\lambda \leq \frac{1}{2L_f\eta_t}$ for all $t\geq 1$, we have
\begin{align}
F(x_{t+1}) - f(x_{t+1},y_{t+1})
& \leq (1-\frac{\eta_t\lambda\mu}{2}) \big(F(x_t) -f(x_t,y_t)\big) + \frac{\eta_t}{8\gamma}\|\tilde{x}_{t+1}-x_t\|^2  -\frac{\eta_t}{4\lambda}\|\tilde{y}_{t+1}-y_t\|^2 \nonumber \\
& \quad + \eta_t\lambda\|\nabla_y f(x_t,y_t)-v_t\|^2,
\end{align}
where $F(x_t)=f(x_t,y^*(x_t))$ with $y^*(x_t) \in \arg\max_{y}f(x_t,y)$ for all $t\geq 1$.
\end{lemma}

\begin{proof}
Using $L_f$-smoothness of $f(x,\cdot)$, such that
\begin{align}
    f(x_{t+1},y_t) + \langle \nabla_y f(x_{t+1},y_t), y_{t+1}-y_t \rangle - \frac{L_f}{2}\|y_{t+1}-y_t\|^2 \leq f(x_{t+1},y_{t+1}),
\end{align}
then we have
\begin{align} \label{eq:B51}
    f(x_{t+1},y_t) & \leq f(x_{t+1},y_{t+1}) - \langle \nabla_y f(x_{t+1},y_t), y_{t+1}-y_t \rangle + \frac{L_f}{2}\|y_{t+1}-y_t\|^2 \nonumber \\
    & = f(x_{t+1},y_{t+1}) - \eta_t\langle \nabla_y f(x_{t+1},y_t), \tilde{y}_{t+1}-y_t \rangle + \frac{L_f\eta^2_t}{2}\|\tilde{y}_{t+1}-y_t\|^2.
\end{align}

Next, we bound the inner product in \eqref{eq:B51}. According to the line 4 of Algorithm~\ref{alg:1}, i.e.,
$\tilde{y}_{t+1} = y_t + \lambda v_t$, we have
\begin{align} \label{eq:B52}
    & - \eta_t\langle \nabla_y f(x_{t+1},y_t), \tilde{y}_{t+1}-y_t \rangle \nonumber \\
    & = - \eta_t\lambda\langle \nabla_y f(x_{t+1},y_t), v_t \rangle \nonumber \\
    & = -\frac{\eta_t\lambda}{2}\Big( \|\nabla_y f(x_{t+1},y_t)\|^2 + \|v_t\|^2 - \|\nabla_y f(x_{t+1},y_t)-\nabla_y f(x_t,y_t) + \nabla_y f(x_t,y_t)-v_t\|^2 \Big) \nonumber \\
    & \leq -\frac{\eta_t\lambda}{2} \|\nabla_y f(x_{t+1},y_t)\|^2 -\frac{\eta_t}{2\lambda} \|\tilde{y}_{t+1}-y_t\|^2 + \eta_t\lambda L^2_f \|x_{t+1}-x_t\|^2 + \eta_t\lambda\|\nabla_y f(x_t,y_t)-v_t\|^2 \nonumber \\
    & \leq -\eta_t\lambda\mu\big(F(x_{t+1})-f(x_{t+1},y_t)\big)-\frac{\eta_t}{2\lambda} \|\tilde{y}_{t+1}-y_t\|^2 + \eta_t\lambda L^2_f \|x_{t+1}-x_t\|^2 + \eta_t\lambda\|\nabla_y f(x_t,y_t)-v_t\|^2,
\end{align}
where the last inequality is due to the quadratic growth condition of $\mu$-PL functions, i.e.,
\begin{align}
    \|\nabla_y f(x_{t+1},y_t)\|^2 \geq 2\mu\big( \max_{y'}f(x_{t+1},y')-f(x_{t+1},y_t)\big) = 2\mu\big( F(x_{t+1})-f(x_{t+1},y_t)\big).
\end{align}
Substituting \eqref{eq:B52} in \eqref{eq:B51}, we have
\begin{align} \label{eq:B53}
    f(x_{t+1},y_t)
    & \leq f(x_{t+1},y_{t+1})-\eta_t\lambda\mu\big(F(x_{t+1})-f(x_{t+1},y_t)\big)-\frac{\eta_t}{2\lambda} \|\tilde{y}_{t+1}-y_t\|^2 + \eta_t\lambda L^2_f \|x_{t+1}-x_t\|^2 \nonumber \\
    & \quad + \eta_t\lambda\|\nabla_y f(x_t,y_t)-v_t\|^2 + \frac{L_f\eta^2_t}{2}\|\tilde{y}_{t+1}-y_t\|^2,
\end{align}
then rearranging the terms, we can obtain
\begin{align} \label{eq:B54}
    F(x_{t+1}) - f(x_{t+1},y_{t+1})
    & \leq (1-\eta_t\lambda\mu)\big(F(x_{t+1})-f(x_{t+1},y_t)\big)-\frac{\eta_t}{2\lambda} \|\tilde{y}_{t+1}-y_t\|^2 + \eta_t\lambda L^2_f \|x_{t+1}-x_t\|^2 \nonumber \\
    & \quad + \eta_t\lambda\|\nabla_y f(x_t,y_t)-v_t\|^2 + \frac{L_f\eta^2_t}{2}\|\tilde{y}_{t+1}-y_t\|^2.
\end{align}

Next, using $L_f$-smoothness of function $f(\cdot,y_t)$, such that
\begin{align}
    f(x_t,y_t) + \langle \nabla_x f(x_t,y_t), x_{t+1}-x_t \rangle - \frac{L_f}{2}\|x_{t+1}-x_t\|^2 \leq f(x_{t+1},y_t),
\end{align}
then we have
\begin{align}
    & f(x_t,y_t) -f(x_{t+1},y_t)  \nonumber \\
    & \leq -\langle \nabla_x f(x_t,y_t), x_{t+1}-x_t \rangle + \frac{L_f}{2}\|x_{t+1}-x_t\|^2 \nonumber \\
    & = -\eta_t\langle \nabla_x f(x_t,y_t) - \nabla F(x_t), \tilde{x}_{t+1}-x_t \rangle - \eta_t\langle \nabla F(x_t), \tilde{x}_{t+1}-x_t \rangle + \frac{L_f\eta^2_t}{2}\|\tilde{x}_{t+1}-x_t\|^2 \nonumber \\
    & \leq \frac{\eta_t}{8\gamma}\|\tilde{x}_{t+1}-x_t\|^2 + 2\eta_t\gamma\|\nabla_x f(x_t,y_t) - \nabla F(x_t)\|^2  - \eta_t\langle \nabla F(x_t), \tilde{x}_{t+1}-x_t \rangle + \frac{L_f\eta^2_t}{2}\|\tilde{x}_{t+1}-x_t\|^2 \nonumber \\
    & \leq \frac{\eta_t}{8\gamma}\|\tilde{x}_{t+1}-x_t\|^2 + 2L^2_f\eta_t\gamma \|y_t - y^*(x_t)\|^2 + F(x_t) - F(x_{t+1}) \nonumber \\
    &  + \frac{\eta^2_tL}{2}\|\tilde{x}_{t+1}-x_t\|^2 + \frac{\eta^2_tL_f}{2}\|\tilde{x}_{t+1}-x_t\|^2 \nonumber \\
    & \leq \frac{4L^2_f\eta_t\gamma}{\mu} \big(F(x_t) -f(x_t,y_t)\big) + F(x_t) - F(x_{t+1}) + \eta_t(\frac{1}{8\gamma}+\eta_tL)\|\tilde{x}_{t+1}-x_t\|^2,
\end{align}
where the second last inequality is due to Lemma~\ref{lem:A1}, i.e.,
$L$-smoothness of function $F(x)$, and the last inequality holds by Lemma~\ref{lem:A2} and $L_f\leq L$.
Then we have
\begin{align} \label{eq:B55}
    F(x_{t+1}) - f(x_{t+1},y_t) & = F(x_{t+1}) - F(x_t) + F(x_t)- f(x_t,y_t) + f(x_t,y_t) -f(x_{t+1},y_t) \nonumber \\
    & \leq (1+\frac{4L^2_f\eta_t\gamma}{\mu}) \big(F(x_t) -f(x_t,y_t)\big) + \eta_t(\frac{1}{8\gamma}+\eta_tL)\|\tilde{x}_{t+1}-x_t\|^2.
\end{align}

Substituting \eqref{eq:B55} in \eqref{eq:B54}, we get
\begin{align}
    & F(x_{t+1}) - f(x_{t+1},y_{t+1}) \nonumber \\
    & \leq (1-\eta_t\lambda\mu)(1+\frac{4L^2_f\eta_t\gamma}{\mu}) \big(F(x_t) -f(x_t,y_t)\big) + \eta_t(\frac{1}{8\gamma}+\eta_tL)(1-\eta_t\lambda\mu)\|\tilde{x}_{t+1}-x_t\|^2 \nonumber \\
    & \quad -\frac{\eta_t}{2\lambda} \|\tilde{y}_{t+1}-y_t\|^2 + \eta_t\lambda L^2_f \|x_{t+1}-x_t\|^2 + \eta_t\lambda\|\nabla_y f(x_t,y_t)-v_t\|^2 + \frac{L_f\eta^2_t}{2}\|\tilde{y}_{t+1}-y_t\|^2 \nonumber \\
    & = (1-\eta_t\lambda\mu)(1+\frac{4L^2_f\eta_t\gamma}{\mu}) \big(F(x_t) -f(x_t,y_t)\big) + \eta_t\big(\frac{1}{8\gamma}+\eta_tL-\frac{\eta_t\lambda\mu}{8\gamma}-\eta^2_tL\lambda\mu+\eta^2_tL^2_f\lambda\big)\|\tilde{x}_{t+1}-x_t\|^2 \nonumber \\
    & \quad -\frac{\eta_t}{2}\big(\frac{1}{\lambda}-L_f\eta_t\big) \|\tilde{y}_{t+1}-y_t\|^2 + \eta_t\lambda\|\nabla_y f(x_t,y_t)-v_t\|^2 \nonumber \\
    & \leq (1-\frac{\eta_t\lambda\mu}{2}) \big(F(x_t) -f(x_t,y_t)\big) + \frac{\eta_t}{8\gamma}\|\tilde{x}_{t+1}-x_t\|^2  -\frac{\eta_t}{4\lambda}\|\tilde{y}_{t+1}-y_t\|^2 + \eta_t\lambda\|\nabla_y f(x_t,y_t)-v_t\|^2,
\end{align}
where the last inequality holds by $L=L_f(1+\frac{\kappa}{2})$, $\gamma\leq \big(\frac{\lambda\mu}{16L},\frac{\mu}{16L^2_f}\big)$ and $\lambda\leq \frac{1}{2L_f\eta_t}$ for all $t\geq 1$, i.e.,
\begin{align}
   & \gamma\leq \frac{\lambda\mu}{16L} \Rightarrow \lambda \geq \frac{16L\gamma}{\mu} =  16(\kappa+\frac{\kappa^2}{2})\gamma\geq 8\kappa^2\gamma \Rightarrow  \frac{\eta_t\lambda\mu}{2} \geq \frac{4L^2_f\eta_t\gamma}{\mu} \nonumber \\
   & \gamma\leq \big(\frac{\lambda\mu}{16L},\frac{\mu}{16L^2_f}\big), \eta_t\in (0,1)\Rightarrow \frac{\eta_t\lambda\mu}{8\gamma} \geq \eta_tL+\eta^2_tL^2_f\lambda \nonumber \\
   &\lambda \leq \frac{1}{2\eta_tL_f}  \Rightarrow \frac{1}{2\lambda} \geq \eta_t L_f, \ \forall t\geq 1.
\end{align}

\end{proof}

\begin{lemma} \label{lem:D1}
 Suppose that the sequence $\big\{x_t,\tilde{x}_t\big\}_{t=1}^T$ be generated from Algorithm \ref{alg:1}.
 Let $0<\eta_t \leq 1$ and $0< \gamma \leq \frac{1}{2L\eta_t}$,
 then we have
 \begin{align}
  \Psi(x_{t+1}) \leq \Psi(x_t) - \frac{\eta_t\gamma}{2}\|\mathcal{G}(x_t,w_t,\gamma)\|^2 +\frac{4\eta_t\gamma L^2_f}{\mu}\big(F(x_t)-f(x_t,y_t)\big) + 2\eta_t\gamma\|\nabla_x f(x_t,x_t)-w_t\|^2,
 \end{align}
where $\mathcal{G}(x_t,w_t,\gamma) = \frac{1}{\gamma}(x_t-\tilde{x}_{t+1})$ and $\Psi(x) = F(x)+\psi(x)$.
\end{lemma}

\begin{proof}
By the line 4 of Algorithm~\ref{alg:1}, we have
\begin{align} \label{eq:D1}
\tilde{x}_{t+1} = \mathcal{P}^{\gamma}_{\psi(\cdot)}\big(x_t-\gamma w_t\big) = \arg\min_{x\in \mathbb{R}^d}\Big\{ \langle w_t, x\rangle + \frac{1}{2\gamma}\|x-x_t\|^2 + \psi(x)\Big\}.
\end{align}

By the optimality condition of the subproblem~(\ref{eq:D1}), we have for any $z\in \mathbb{R}^d$
\begin{align}
 \big\langle w_t + \frac{1}{\gamma}(\tilde{x}_{t+1}-x_t) + \nu_{t+1}, z-\tilde{x}_{t+1}\big\rangle \geq 0,
\end{align}
where $\nu_{t+1}\in \partial \psi(\tilde{x}_{t+1})$.

By using the convexity of $\psi(x)$, and
let $z=x_t$, we can obtain
\begin{align}  \label{eq:D2}
 \langle w_t, x_t - \tilde{x}_{t+1}\rangle & \geq \frac{1}{\gamma}\|\tilde{x}_{t+1}-x_t\|^2 + \langle \nu_{t+1}, \tilde{x}_{t+1}-x_t\rangle \nonumber \\
 & \geq \frac{1}{\gamma}\|\tilde{x}_{t+1}-x_t\|^2 + \psi(\tilde{x}_{t+1})-\psi(x_t).
\end{align}

Let $\mathcal{G}(x_t,w_t,\gamma) = \frac{1}{\gamma}(x_t-\tilde{x}_{t+1})$. According to Lemma~\ref{lem:A1}, i.e., function $F(x)$ is $L$-smooth, we have
\begin{align} \label{eq:D3}
  F(x_{t+1}) & \leq F(x_t) + \langle \nabla F(x_t), x_{t+1}-x_t\rangle + \frac{L}{2}\|x_{t+1}-x_t\|^2 \nonumber \\
  & = F(x_t) + \eta_t\langle w_t, \tilde{x}_{t+1}-x_t\rangle + \eta_t\gamma \langle \nabla F(x_t) -w_t, \mathcal{G}(x_t,w_t,\gamma)\rangle+ \frac{\gamma^2\eta_t^2L}{2}\|\mathcal{G}(x_t,w_t,\gamma)\|^2 \nonumber \\
  & \mathop{\leq}^{(i)} F(x_t) - \gamma\eta_t \|\mathcal{G}(x_t,w_t,\gamma)\|^2 - \eta_t\psi(\tilde{x}_{t+1}) + \eta_t\psi(x_t) + \eta_t\gamma \langle \nabla F(x_t) - w_t, \mathcal{G}(x_t,w_t,\gamma)\rangle \nonumber \\
  & \quad + \frac{\gamma ^2\eta_t^2L}{2}\|\mathcal{G}(x_t,w_t,\gamma)\|^2 \nonumber \\
  & = F(x_t) - \gamma\eta_t \|\mathcal{G}(x_t,w_t,\gamma)\|^2 - \eta_t\psi(\tilde{x}_{t+1}) - (1-\eta_t)\psi(x_t) + \psi(x_t) \nonumber \\
  & \quad  + \eta_t\gamma \langle \nabla F(x_t) - w_t, \mathcal{G}(x_t,w_t,\gamma)\rangle+ \frac{\gamma ^2\eta_t^2L}{2}\|\mathcal{G}(x_t,w_t,\gamma)\|^2 \nonumber \\
  & \mathop{\leq}^{(ii)} F(x_t) - \gamma\eta_t \|\mathcal{G}(x_t,w_t,\gamma)\|^2 - \psi(x_{t+1}) + \psi(x_t) \nonumber \\
  & \quad  + \eta_t\gamma \langle \nabla F(x_t) - w_t, \mathcal{G}(x_t,w_t,\gamma)\rangle+ \frac{\gamma ^2\eta_t^2L}{2}\|\mathcal{G}(x_t,w_t,\gamma)\|^2 \nonumber \\
  & \mathop{\leq}^{(iii)} F(x_t) - \frac{\eta_t\gamma}{2}\|\mathcal{G}(x_t,w_t,\gamma)\|^2 - \psi(x_{t+1}) + \psi(x_t) + \eta_t\gamma\|w_t - \nabla F(x_t)\|^2,
\end{align}
where the inequality (i) holds by the above inequality~(\ref{eq:D2}), and the inequality (ii) is due to $x_{t+1} = x_t + \eta_t(\tilde{x}_{t+1}-x_t)$ and the convexity of function $\psi(x)$, i.e., $\psi(x_{t+1})=\psi((1-\eta_t)x_t + \eta_t\tilde{x}_{t+1})\leq (1-\eta_t)\psi(x_t) + \eta_t\psi(\tilde{x}_{t+1})$,
and the last inequality (iii) holds by $0<\gamma\leq \frac{1}{2\eta_t L}$ and the following inequality
\begin{align}
 \langle \nabla F(x_t) -w_t, \mathcal{G}(x_t,w_t,\gamma)\rangle &\leq \|w_t - \nabla F(x_t)\|\|\mathcal{G}(x_t,w_t,\gamma)\| \nonumber \\
 & \leq  \|w_t - \nabla F(x_t)\|^2 + \frac{1}{4}\|\mathcal{G}(x_t,w_t,\gamma)\|^2,
\end{align}
where the above inequality holds by Young inequality.

Considering the bound of the term $\|\nabla F(x_t)-w_t\|^2$, then we have
\begin{align} \label{eq:D4}
  \|\nabla F(x_t)-w_t\|^2 & = \|\nabla_x f(x_t,y^*(x_t))- \nabla_x f(x_t,y_t)+\nabla_x f(x_t,y_t)- w_t\|^2\nonumber \\
  & \leq 2\|\nabla_x f(x_t,y^*(x_t))- \nabla_x f(x_t,y_t)\|^2 + 2\|\nabla_x f(x_t,x_t)-w_t\|^2
  \nonumber \\
  & \leq 2L^2_f\|y^*(x_t) - y_t\|^2 + 2\|\nabla_x f(x_t,x_t)-w_t\|^2,
\end{align}
where the first inequality is due to the Cauchy-Schwarz inequality and the second is due to Young's inequality.

Let $\Psi(x) = F(x)+\psi(x)$.
By plugging the above inequalities \eqref{eq:D4} into \eqref{eq:D3},
we obtain
\begin{align} \label{eq:D5}
  \Psi(x_{t+1}) &\leq \Psi(x_t) - \frac{\eta_t\gamma}{2}\|\mathcal{G}(x_t,w_t,\gamma)\|^2 + \eta_t\gamma\|w_t - \nabla F(x_t)\|^2 \\
  & \leq  \Psi(x_t) - \frac{\eta_t\gamma}{2}\|\mathcal{G}(x_t,w_t,\gamma)\|^2 +2\eta_t\gamma L^2_f\|y^*(x_t) - y_t\|^2 + 2\eta_t\gamma\|\nabla_x f(x_t,x_t)-w_t\|^2   \nonumber \\
  & \leq \Psi(x_t) - \frac{\eta_t\gamma}{2}\|\mathcal{G}(x_t,w_t,\gamma)\|^2 +\frac{4\eta_t\gamma L^2_f}{\mu}\big(F(x_t)-f(x_t,y_t)\big) + 2\eta_t\gamma\|\nabla_x f(x_t,x_t)-w_t\|^2, \nonumber
\end{align}
where the last inequality holds by the above Lemma~\ref{lem:A2} using in $F(x_t)=f(x_t,y^*(x_t))=\max_y f(x_t,y)$ with $y^*(x_t)\in \arg\max f(x_t,y)$.

\end{proof}

\begin{theorem}  \label{th:A1}
(Restatement of Theorem 1)
 Under the above Assumptions~\ref{ass:1}-\ref{ass:4}, in Algorithm~\ref{alg:1}, let $\eta_t=\frac{k}{(m+t)^{1/3}}$ for all $t\geq 0$, $\alpha_{t+1}=c_1\eta_t^2$, $\beta_{t+1}=c_2\eta_t^2$, $m \geq \max\big(2,k^3, (c_1k)^3, (c_2k)^3\big)$, $k>0$, $ c_1 \geq \frac{2}{3k^3} + \frac{9L^2_f}{\mu^2}$, $c_2 \geq \frac{2}{3k^3} + \frac{9}{4}$,
 $0< \lambda \leq \min\big(\frac{3}{4\sqrt{2}\mu},\frac{m^{1/3}}{2L_fk}\big)$ and $0< \gamma \leq \big(\frac{\lambda\mu}{16L},\frac{\mu}{16L^2_f},\frac{m^{1/3}}{2Lk},\frac{1}{8L_f},\frac{\lambda\mu^2}{9L^2_f}\big)$, we have
\begin{align}
  \frac{1}{T}\sum_{t=1}^T\mathbb{E}\|\mathcal{G}(x_t,\nabla F(x_t),\gamma)\| \leq \frac{2\sqrt{3H}m^{1/6}}{T^{1/2}} + \frac{2\sqrt{3H}}{T^{1/3}},
\end{align}
where $\Psi(x)=F(x)+\psi(x)$ and $H = \frac{\Psi(x_1) - \Psi^*}{\gamma k} + \frac{9 L^2_f}{k\lambda\mu^2}\Delta_1 + \frac{2\sigma^2m^{1/3}}{k^2} + 2k^2(c_1^2+c_2^2)\sigma^2\ln(m+T)$ with $\Delta_1=F(x_1)-f(x_1,y_1)$.
\end{theorem}

\begin{proof}
Since $\eta_t$ is decreasing and $m\geq k^3$, we have $\eta_t \leq \eta_0 = \frac{k}{m^{1/3}} \leq 1$ and $\gamma \leq \frac{1}{2L\eta_0} =\frac{m^{1/3}}{2Lk} \leq \frac{1}{2L\eta_t}$ for any $t\geq 0$.
Meanwhile, we have $\lambda \leq \frac{m^{1/3}}{2L_fk}=\frac{1}{2L_f\eta_0} \leq \frac{1}{2L_f\eta_t}$ for any $t\geq 0$.
Due to $0 < \eta_t \leq 1$ and $m\geq \max\big( (c_1k)^3, (c_2k)^3 \big)$, we have $\alpha_{t+1} = c_1\eta_t^2 \leq c_1\eta_t \leq \frac{c_1k}{m^{1/3}}\leq 1$ and $\beta_{t+1} = c_2\eta_t^2 \leq c_2\eta_t \leq \frac{c_2k}{m^{1/3}}\leq 1$.

According to the above Lemma~\ref{lem:B1}, we can obtain
\begin{align}
& \frac{1}{\eta_t}\mathbb{E} \|\nabla_{x} f(x_{t+1},y_{t+1}) - w_{t+1}\|^2 - \frac{1}{\eta_{t-1}}\mathbb{E} \|\nabla_{x} f(x_t,y_t) - w_t\|^2  \\
& \leq \big(\frac{1 - \beta_{t+1}}{\eta_t}  -  \frac{1}{\eta_{t-1}}\big)\mathbb{E} \|\nabla_{x} f(x_t,y_t) -w_t\|^2 + 4L^2_f\eta_t\mathbb{E}\big(\|\tilde{x}_{t+1}-x_t\|^2 + \|\tilde{y}_{t+1}-y_t\|^2\big) +  \frac{2\beta_{t+1}^2\sigma^2}{\eta_t} \nonumber \\
& = \big(\frac{1}{\eta_t} -  \frac{1}{\eta_{t-1}} - c_2\eta_t\big)\mathbb{E} \|\nabla_{x} f(x_t,y_t) -w_t\|^2 + 4L^2_f\eta_t\mathbb{E}\big(\|\tilde{x}_{t+1}-x_t\|^2 + \|\tilde{y}_{t+1}-y_t\|^2\big) +  2c^2_2\eta_t^3\sigma^2, \nonumber
\end{align}
where the second inequality is due to $0<\beta_{t+1}\leq 1$ and $\beta_{t+1}=c_2\eta_t^2$.
Similarly, we have
\begin{align}
& \frac{1}{\eta_t}\mathbb{E} \|\nabla_{y} f(x_{t+1},y_{t+1}) - v_{t+1}\|^2
-  \frac{1}{\eta_{t-1}}\mathbb{E} \|\nabla_{y} f(x_t,y_t) - v_t\|^2  \\
& \leq \big(\frac{1}{\eta_t} - \frac{1}{\eta_{t-1}} - c_1\eta_t\big)\mathbb{E} \|\nabla_{y} f(x_t,y_t) -v_t\|^2 + 4L^2_f\eta_t\mathbb{E}\big(\|\tilde{x}_{t+1}-x_t\|^2 + \|\tilde{y}_{t+1}-y_t\|^2\big) +
 2c^2_1\eta_t^3\sigma^2. \nonumber
\end{align}
By $\eta_t = \frac{k}{(m+t)^{1/3}}$, we have
\begin{align}
\frac{1}{\eta_t} - \frac{1}{\eta_{t-1}} &= \frac{1}{k}\big( (m+t)^{\frac{1}{3}} - (m+t-1)^{\frac{1}{3}}\big) \nonumber \\
& \leq \frac{1}{3k(m+t-1)^{2/3}} = \frac{2^{2/3}}{3k\big(2(m+t-1)\big)^{2/3}} \nonumber \\
& \leq \frac{2^{2/3}}{3k(m+t)^{2/3}} = \frac{2^{2/3}}{3k^3}\frac{k^2}{(m+t)^{2/3}}= \frac{2^{2/3}}{3k^3}\eta_t^2 \leq \frac{2}{3k^3}\eta_t,
\end{align}
where the first inequality holds by the concavity of function $f(x)=x^{1/3}$, \emph{i.e.}, $(x+y)^{1/3}\leq x^{1/3} + \frac{y}{3x^{2/3}}$, and
the last inequality is due to $0<\eta_t\leq 1$.

Let $c_2 \geq \frac{2}{3k^3} + \frac{9}{4}$, we have
\begin{align} \label{eq:G21}
& \frac{1}{\eta_t}\mathbb{E} \|\nabla_{x} f(x_{t+1},y_{t+1}) - w_{t+1}\|^2 - \frac{1}{\eta_{t-1}}\mathbb{E} \|\nabla_{x} f(x_t,y_t) - w_t\|^2  \\
& \leq -\frac{9\eta_t}{4}\mathbb{E} \|\nabla_{x} f(x_t,y_t) -w_t\|^2 + 4L^2_f\eta_t\mathbb{E}\big(\|\tilde{x}_{t+1}-x_t\|^2 + \|\tilde{y}_{t+1}-y_t\|^2\big) + 2c^2_2\eta_t^3\sigma^2. \nonumber
\end{align}
Let $c_1 \geq \frac{2}{3k^3} + \frac{9L^2_f}{\mu^2}$, we have
\begin{align} \label{eq:G22}
& \frac{1}{\eta_t}\mathbb{E} \|\nabla_{y} f(x_{t+1},y_{t+1}) - w_{t+1}\|^2
-  \frac{1}{\eta_{t-1}}\mathbb{E} \|\nabla_{y} f(x_t,y_t) - w_t\|^2  \\
& \leq - \frac{9L^2_f\eta_t}{\mu^2}\mathbb{E} \|\nabla_{y} f(x_t,y_t) -w_t\|^2 + 4L^2_f\eta_t\mathbb{E}\big(\|\tilde{x}_{t+1}-x_t\|^2 + \|\tilde{y}_{t+1}-y_t\|^2\big) +
2c^2_1\eta_t^3\sigma^2. \nonumber
 \end{align}

According to Lemma \ref{lem:C1}, we have
\begin{align} \label{eq:G23}
F(x_{t+1}) -f(x_{t+1},y_{t+1}) - \big(F(x_t) -f(x_t,y_t)\big)  & \leq -\frac{\eta_t\lambda\mu}{2} \big(F(x_t) -f(x_t,y_t)\big) + \frac{\eta_t}{8\gamma}\|\tilde{x}_{t+1}-x_t\|^2  \nonumber \\
& \quad -\frac{\eta_t}{4\lambda}\|\tilde{y}_{t+1}-y_t\|^2 + \eta_t\lambda\|\nabla_y f(x_t,y_t)-v_t\|^2.
\end{align}
According to Lemma \ref{lem:D1}, we have
\begin{align} \label{eq:G24}
\Psi(x_{t+1})-\Psi(x_t) \leq - \frac{\eta_t\gamma}{2}\|\mathcal{G}(x_t,w_t,\gamma)\|^2 +\frac{4\eta_t\gamma L^2_f}{\mu}\big(F(x_t)-f(x_t,y_t)\big) + 2\eta_t\gamma\|\nabla_x f(x_t,x_t)-w_t\|^2.
\end{align}

Next, we define a Lyapunov function (i.e., potential function), for any $t\geq 1$
\begin{align}
\Omega_t = \mathbb{E}\big[ \Psi(x_t) + \frac{9\gamma L^2_f}{\lambda\mu^2}\big(F(x_t) -f(x_t,y_t) \big) + \gamma\big(\frac{1}{\eta_{t-1}}\|\nabla_{x} f(x_t,y_t)-w_t\|^2 + \frac{1}{\eta_{t-1}}\|\nabla_{y} f(x_t,y_t)-v_t\|^2 \big) \big].
\end{align}
Let $\mathcal{G}(x_t,w_t,\gamma) = \frac{1}{\gamma}(x_t-\tilde{x}_{t+1})$, then we have
\begin{align}
& \Omega_{t+1}- \Omega_t \nonumber \\
& = \Psi(x_{t+1})  -  \Psi(x_t) + \frac{9\gamma L^2_f}{\lambda\mu^2}
\big(F(x_{t+1}) -f(x_{t+1},y_{t+1}) - (F(x_t) -f(x_t,y_t)) \big)
 + \gamma \big(\frac{1}{\eta_t}\mathbb{E}\|\nabla_{x} f(x_{t+1},y_{t+1})-w_{t+1}\|^2 \nonumber \\
& \quad  -  \frac{1}{\eta_{t-1}}\mathbb{E}\|\nabla_{x} f(x_t,y_t)-w_t\|^2  +  \frac{1}{\eta_t}\mathbb{E}\|\nabla_{y} f(x_{t+1},y_{t+1})-v_{t+1}\|^2
 -  \frac{1}{\eta_{t-1}}\mathbb{E}\|\nabla_{y} f(x_t,y_t)-v_t\|^2 \big) \nonumber \\
& \leq \frac{4\gamma L^2_f\eta_t}{\mu}\big(F(x_t)-f(x_t,y_t)\big) + 2\gamma\eta_t\|\nabla_x f(x_t,x_t)-w_t\|^2 - \frac{\eta_t\gamma}{2}\|\mathcal{G}(x_t,w_t,\gamma)\|^2 \nonumber \\
& \quad  + \frac{9\gamma L^2_f}{\lambda\mu^2} \bigg( -\frac{\eta_t\lambda\mu}{2} \big(F(x_t) -f(x_t,y_t)\big) + \frac{\eta_t\gamma}{8}\|\mathcal{G}(x_t,w_t,\gamma)\|^2  -\frac{\eta_t}{4\lambda}\|\tilde{y}_{t+1}-y_t\|^2 + \eta_t\lambda\|\nabla_y f(x_t,y_t)-v_t\|^2 \bigg)  \nonumber \\
& \quad  -\frac{9\gamma\eta_t}{4}\mathbb{E} \|\nabla_{x} f(x_t,y_t) -w_t\|^2 + 4\gamma L^2_f\eta_t \mathbb{E}\big(\gamma^2\|\mathcal{G}(x_t,w_t,\gamma)\|^2 + \|\tilde{y}_{t+1}-y_t\|^2\big) + 2\gamma c^2_2\eta_t^3\sigma^2\nonumber \\
& \quad -\frac{9\gamma L^2_f\eta_t}{\mu^2} \mathbb{E}\|\nabla_{y} f(x_t,y_t) -v_t\|^2  + 4\gamma L^2_f\eta_t\mathbb{E}\big(\gamma^2\|\mathcal{G}(x_t,w_t,\gamma)\|^2 + \|\tilde{y}_{t+1}-y_t\|^2\big) + 2\gamma c^2_1\eta_t^3\sigma^2 \nonumber \\
& \leq -\frac{\gamma L_f^2\eta_t}{2\mu}\big(F(x_t) -f(x_t,y_t) \big) - \frac{\gamma\eta_t}{4}\mathbb{E}\|\nabla_{x} f(x_t,y_t) -w_t\|^2 +2(c_1^2+c^2_2)\gamma \sigma^2\eta_t^3 \nonumber \\
& \quad -\big(\frac{9\gamma L^2_f}{4\lambda^2\mu^2} - 8\gamma L^2_f\big)\eta_t\mathbb{E}\|\tilde{y}_{t+1} - y_t\|^2 - \big( \frac{\gamma}{2} - 8\gamma^3 L^2_f - \frac{9\gamma^2 L^2_f}{8\lambda\mu^2}\big)\eta_t\mathbb{E}\|\mathcal{G}(x_t,w_t,\gamma)\|^2 \nonumber \\
& \leq -\frac{\gamma L_f^2\eta_t}{2\mu}\big(F(x_t) -f(x_t,y_t) \big)  -  \frac{\gamma\eta_t}{4}\mathbb{E}\|\nabla_{x} f(x_t,y_t) -w_t\|^2 - \frac{\eta_t\gamma}{4}\mathbb{E}\|\mathcal{G}(x_t,w_t,\gamma)\|^2 +2(c_1^2+c^2_2)\gamma \sigma^2\eta_t^3,
\end{align}
where the first inequality holds by the above inequalities \eqref{eq:G21}, \eqref{eq:G22}, \eqref{eq:G23} and \eqref{eq:G24}; the last inequality is due to $0< \lambda \leq \frac{3}{4\sqrt{2}\mu}$ and $0< \gamma \leq \big(\frac{1}{8L_f},\frac{\lambda\mu^2}{9L^2_f}\big)$.
Thus, we have
\begin{align} \label{eq:G25}
 & \frac{L_f^2\eta_t}{2\mu}\big(F(x_t) -f(x_t,y_t) \big) + \frac{\eta_t}{4}\mathbb{E}\|\nabla_{x} f(x_t,y_t) -w_t\|^2 + \frac{\eta_t}{4}\mathbb{E}\|\mathcal{G}(x_t,w_t,\gamma)\|^2 \nonumber \\
 & \leq \frac{\Omega_t - \Omega_{t+1}}{\gamma} + 2(c_1^2+c^2_2)\sigma^2\eta_t^3.
\end{align}

Taking average over $t=1,2,\cdots,T$ on both sides of \eqref{eq:G25}, we have
\begin{align}
 & \frac{1}{T} \sum_{t=1}^T \Big( \frac{L_f^2\eta_t}{2\mu}\big(F(x_t) -f(x_t,y_t) \big) + \frac{\eta_t}{4}\mathbb{E}\|\nabla_{x} f(x_t,y_t) -w_t\|^2 + \frac{\eta_t}{4}\mathbb{E}\|\mathcal{G}(x_t,w_t,\gamma)\|^2  \Big) \nonumber \\
 & \leq  \sum_{t=1}^T \frac{\Omega_t - \Omega_{t+1}}{T\gamma} + \frac{2(c_1^2+c^2_2) \sigma^2}{T}\sum_{t=1}^T\eta_t^3. \nonumber
\end{align}
Let $\Delta_1 = F(x_1) -f(x_1,y_1)$, we have
\begin{align} \label{eq:G26}
 \Omega_1 &= \Psi(x_1) + \frac{9\gamma L^2_f}{\lambda\mu^2}\big( F(x_1) -f(x_1,y_1)\big) + \gamma\big(\frac{1}{\eta_0}\|\nabla_{x}f(x_1,y_1)-w_1\|^2 + \frac{1}{\eta_0}\|\nabla_{y} f(x_1,y_1)-v_1\|^2 \big)  \nonumber \\
 & \leq \Psi(x_1) + \frac{9\gamma L^2_f}{\lambda\mu^2}\Delta_1  + \frac{2\gamma\sigma^2}{\eta_0},
\end{align}
where the last inequality holds by Assumption~\ref{ass:3}.

Since $\eta_t$ is decreasing, i.e., $\eta_T^{-1} \geq \eta_t^{-1}$ for any $0\leq t\leq T$, we have
 \begin{align}
 & \frac{1}{T} \sum_{t=1}^T \mathbb{E}\Big[  \frac{L_f^2}{2\mu}\big(F(x_t) -f(x_t,y_t) \big) + \frac{1}{4}\|\nabla_{x} f(x_t,y_t) -w_t\|^2 + \frac{1}{4}\|\mathcal{G}(x_t,w_t,\gamma)\|^2 \Big] \nonumber \\
 & \leq  \sum_{t=1}^T \frac{\Omega_t - \Omega_{t+1}}{\eta_TT\gamma} + \frac{2(c_1^2+c^2_2) \sigma^2}{\eta_TT}\sum_{t=1}^T\eta_t^3 \nonumber \\
 & \leq \frac{1}{\eta_TT\gamma}\Big(\Psi(x_1) + \frac{9\gamma L^2_f}{\lambda\mu^2}\Delta_1  + \frac{2\gamma\sigma^2}{\eta_0} - \Psi^* \Big) + \frac{2(c_1^2+c^2_2) \sigma^2}{\eta_TT}\int^T_1\frac{k^3}{m+t}dt \nonumber \\
 & \leq \frac{\Psi(x_1)-\Psi^*}{\eta_TT\gamma} + \frac{9 L^2_f}{\eta_TT\lambda\mu^2}\Delta_1  + \frac{2\sigma^2}{\eta_TT\eta_0} + \frac{2k^3(c_1^2+c^2_2) \sigma^2}{\eta_TT}\ln(m+T) \nonumber \\
 & = \bigg( \frac{\Psi(x_1) - \Psi^*}{\gamma k} + \frac{9 L^2_f}{k\lambda\mu^2}\Delta_1  + \frac{2\sigma^2m^{1/3}}{k^2} + 2k^2(c_1^2+c_2^2)\sigma^2\ln(m+T)\bigg) \frac{(m+T)^{1/3}}{T},
\end{align}
where the second inequality holds by the above inequality \eqref{eq:G26}. Let $H = \frac{\Psi(x_1) - \Psi^*}{\gamma k} + \frac{9 L^2_f}{k\lambda\mu^2}\Delta_1 + \frac{2\sigma^2m^{1/3}}{k^2} + 2k^2(c_1^2+c_2^2)\sigma^2\ln(m+T)$,
we have
\begin{align} \label{eq:G27}
 \frac{1}{T} \sum_{t=1}^T \mathbb{E}\big[ \frac{L_f^2}{2\mu}\big(F(x_t) -f(x_t,y_t) \big) + \frac{1}{4}\|\nabla_{x} f(x_t,y_t) -w_t\|^2 + \frac{1}{4}\|\mathcal{G}(x_t,w_t,\gamma)\|^2 \big]  \leq \frac{H}{T}(m+T)^{1/3}.
\end{align}

Next, we define a useful gradient mapping
$\mathcal{G}(x_t,\nabla F(x_t),\gamma)=\frac{1}{\gamma}(x_t-x^+_{t+1})$, where $x^+_{t+1}$
is generated from
\begin{align}
x^+_{t+1} = \arg\min_{x\in \mathbb{R}^d}\Big\{ \langle \nabla F(x_t), x\rangle
+ \frac{1}{2\gamma}\|x-x_t\|^2 + \psi(x)\Big\}  = \mathcal{P}_{\psi(\cdot)}^\gamma(x_t-\gamma\nabla F(x_t)),
\end{align}
where $F(x)=f(x,y^*(x))$ with $y^*(x)\in \arg\max_y f(x,y)$.

According to Lemma~\ref{lem:A1}, we have $\nabla F(x) =\nabla_{x} f(x_t,y^*(x_t))$. Then we
obtain
\begin{align} \label{eq:G28}
 \|\mathcal{G}(x_t,\nabla F(x_t),\gamma)\| & = \|\mathcal{G}(x_t,\nabla F(x_t),\gamma)-\mathcal{G}(x_t,w_t,\gamma)+\mathcal{G}(x_t,w_t,\gamma)\| \nonumber \\
 & \leq \|\mathcal{G}(x_t,\nabla F(x_t),\gamma)-\mathcal{G}(x_t,w_t,\gamma)\| + \|\mathcal{G}(x_t,w_t,\gamma)\| \nonumber \\
 & = \|\frac{1}{\gamma}(x_t-\mathcal{P}_{\psi(\cdot)}^\gamma(x_t-\gamma \nabla F(x_t)))-\frac{1}{\gamma}(x_t-\mathcal{P}_{\psi(\cdot)}^\gamma(x_t-\gamma w_t))\| + \|\mathcal{G}(x_t,w_t,\gamma)\| \nonumber \\
 & \mathop{\leq}^{(i)} \|w_t - \nabla F(x_t)\|+ \|\mathcal{G}(x_t,w_t,\gamma)\| \nonumber \\
 & \leq \|w_t - \nabla_{x} f(x_t,y_t)\| + \|\nabla_{x} f(x_t,y_t) - \nabla_{x} f(x_t,y^*(x_t)) \|+ \|\mathcal{G}(x_t,w_t,\gamma)\| \nonumber \\
 & \leq \|w_t - \nabla_{x} f(x_t,y_t)\| + L_f\|y_t - y^*(x_t) \|+ \|\mathcal{G}(x_t,w_t,\gamma)\| \nonumber \\
 & \leq \|w_t - \nabla_{x} f(x_t,y_t)\| + \frac{\sqrt{2}L_f}{\sqrt{\mu}}\sqrt{F(x_t)-f(x_t,y_t)} + \|\mathcal{G}(x_t,w_t,\gamma)\|,
\end{align}
where the inequality~(i) holds by the lemma 2 of \cite{ghadimi2016mini}, and the last inequality holds by Lemma~\ref{lem:A2}.

According to the above inequalities~(\ref{eq:G27}) and~(\ref{eq:G28}), we have
\begin{align} \label{eq:G29}
 & \frac{1}{T}\sum_{t=1}^T\mathbb{E}\|\mathcal{G}(x_t,\nabla F(x_t),\gamma)\| \nonumber \\
 & \leq \frac{1}{T}\sum_{t=1}^T
\mathbb{E}\big[\|w_t - \nabla_{x} f(x_t,y_t)\| + \frac{\sqrt{2}L_f}{\sqrt{\mu}}\sqrt{F(x_t)-f(x_t,y_t)} + \|\mathcal{G}(x_t,w_t,\gamma)\|\big] \nonumber \\
 & \mathop{\leq}^{(i)} \Big( \frac{1}{T} \sum_{t=1}^T \mathbb{E}\big[ \frac{6L_f^2}{\mu}\big(F(x_t)-f(x_t,y_t)\big) + 3\|\nabla_{x} f(x_t,y_t) -w_t\|^2 +3\|\mathcal{G}(x_t,w_t,\gamma)\|^2 \big] \Big)^{1/2} \nonumber \\
 & \leq \frac{2\sqrt{3H}}{\sqrt{T}}(m+T)^{1/6}\leq \frac{2\sqrt{3H}m^{1/6}}{T^{1/2}} + \frac{2\sqrt{3H}}{T^{1/3}},
\end{align}
where the inequality (i) holds by Jensen's inequality.

\end{proof}

\subsection{ Convergence Analysis of AdaMSGDA Algorithm }
In this subsection, we provide the detailed convergence analysis of AdaMSGDA algorithm
under some mild conditions.

\begin{lemma} \label{lem:C2}
(Restatement of Lemma 3)
Suppose the sequence $\{x_t,y_t\}_{t=1}^T$ be generated from Algorithm~\ref{alg:2}.
Under the above Assumptions, given $\gamma\leq \min\big(\frac{\lambda\mu}{16\rho_uL},\frac{\rho_l\mu}{16\rho_uL^2_f}\big)$ and $\lambda \leq \frac{1}{2\eta_t L_f\rho_u}$ for all $t\geq 1$, we have
\begin{align}
F(x_{t+1}) - f(x_{t+1},y_{t+1})
& \leq (1-\frac{\eta_t\lambda\mu}{2\rho_u}) \big(F(x_t) -f(x_t,y_t)\big) + \frac{\eta_t}{8\gamma}\|\tilde{x}_{t+1}-x_t\|^2  -\frac{\eta_t}{4\lambda\rho_u}\|\tilde{y}_{t+1}-y_t\|^2 \nonumber \\
& \quad + \frac{\eta_t\lambda}{\rho_l}\|\nabla_y f(x_t,y_t)-v_t\|^2,
\end{align}
where $F(x_t)=f(x_t,y^*(x_t))$ with $y^*(x_t) \in \arg\min_{y}f(x_t,y)$ for all $t\geq 1$.
\end{lemma}

\begin{proof}
Using $L_f$-smoothness of $f(x,\cdot)$, such that
\begin{align}
    f(x_{t+1},y_t) + \langle \nabla_y f(x_{t+1},y_t), y_{t+1}-y_t \rangle - \frac{L_f}{2}\|y_{t+1}-y_t\|^2 \leq f(x_{t+1},y_{t+1}),
\end{align}
then we have
\begin{align} \label{eq:B61}
    f(x_{t+1},y_t) & \leq f(x_{t+1},y_{t+1}) - \langle \nabla_y f(x_{t+1},y_t), y_{t+1}-y_t \rangle + \frac{L_f}{2}\|y_{t+1}-y_t\|^2 \nonumber \\
    & = f(x_{t+1},y_{t+1}) - \eta_t\langle \nabla_y f(x_{t+1},y_t), \tilde{y}_{t+1}-y_t \rangle + \frac{L_f\eta^2_t}{2}\|\tilde{y}_{t+1}-y_t\|^2.
\end{align}
Since $\rho_u I_{p} \succeq B_t \succeq \rho_l I_{p} \succ 0$ for any $t\geq 1$ is positive definite, we set $B_t=L_t(L_t)^T$, where $ \sqrt{\rho_u} I_{p} \succeq L_t \succeq \sqrt{\rho_l} I_{p} \succ 0$. Thus, we have $B^{-1}_t=(L^{-1}_t)^{T}L^{-1}_t$, where $ \frac{1}{\sqrt{\rho_l}}I_{p} \succeq L^{-1}_t \succeq \frac{1}{\sqrt{\rho_u}}I_{p} \succ 0$.

Next, we bound the inner product in \eqref{eq:B61}.
According to the line 5 of Algorithm~\ref{alg:2}, i.e.,
$\tilde{y}_{t+1} = y_t + \lambda B_t^{-1}v_t$, we have
\begin{align} \label{eq:B62}
    & - \eta_t\langle \nabla_y f(x_{t+1},y_t), \tilde{y}_{t+1}-y_t \rangle \nonumber \\
    & = - \eta_t\lambda\langle \nabla_y f(x_{t+1},y_t), B^{-1}_tv_t \rangle \nonumber \\
    & = - \eta_t\lambda\langle L^{-1}_t\nabla_y f(x_{t+1},y_t), L^{-1}_tv_t \rangle \nonumber \\
    & = -\frac{\eta_t\lambda}{2}\Big( \|L^{-1}_t\nabla_y f(x_{t+1},y_t)\|^2 + \|L^{-1}_tv_t\|^2 - \|L^{-1}_t\nabla_y f(x_{t+1},y_t)-L^{-1}_t\nabla_y f(x_t,y_t) + L^{-1}_t\nabla_y f(x_t,y_t)-L^{-1}_tv_t\|^2 \Big) \nonumber \\
    & \leq -\frac{\eta_t\lambda}{2\rho_u} \|\nabla_y f(x_{t+1},y_t)\|^2 -\frac{\eta_t}{2\lambda\rho_u} \|\tilde{y}_{t+1}-y_t\|^2 + \frac{\eta_t\lambda L^2_f}{\rho_l} \|x_{t+1}-x_t\|^2 + \frac{\eta_t\lambda}{\rho_l}\|\nabla_y f(x_t,y_t)-v_t\|^2 \nonumber \\
    & \leq -\frac{\eta_t\lambda\mu}{\rho_u}\big(F(x_{t+1})-f(x_{t+1},y_t)\big)-\frac{\eta_t}{2\lambda\rho_u} \|\tilde{y}_{t+1}-y_t\|^2 + \frac{\eta_t\lambda L^2_f}{\rho_l} \|x_{t+1}-x_t\|^2 + \frac{\eta_t\lambda}{\rho_l}\|\nabla_y f(x_t,y_t)-v_t\|^2,
\end{align}
where the last inequality is due to the quadratic growth condition of $\mu$-PL functions, i.e.,
\begin{align}
    \|\nabla_y f(x_{t+1},y_t)\|^2 \geq 2\mu\big( \max_{y'}f(x_{t+1},y')-f(x_{t+1},y_t)\big) = 2\mu\big( F(x_{t+1})-f(x_{t+1},y_t)\big).
\end{align}
Substituting \eqref{eq:B62} into \eqref{eq:B61}, we have
\begin{align} \label{eq:B63}
    f(x_{t+1},y_t)
    & \leq  f(x_{t+1},y_{t+1})-\frac{\eta_t\lambda\mu}{\rho_u}\big(F(x_{t+1})-f(x_{t+1},y_t)\big)-\frac{\eta_t}{2\lambda\rho_u} \|\tilde{y}_{t+1}-y_t\|^2 + \frac{\eta_t\lambda L^2_f}{\rho_l} \|x_{t+1}-x_t\|^2 \nonumber \\
    & \quad + \frac{\eta_t\lambda}{\rho_l}\|\nabla_y f(x_t,y_t)-v_t\|^2 + \frac{L_f\eta^2_t}{2}\|\tilde{y}_{t+1}-y_t\|^2,
\end{align}
then rearranging the terms, we can obtain
\begin{align} \label{eq:B64}
    F(x_{t+1}) - f(x_{t+1},y_{t+1})
    & \leq (1-\frac{\eta_t\lambda\mu}{\rho_u})\big(F(x_{t+1})-f(x_{t+1},y_t)\big)-\frac{\eta_t}{2\lambda\rho_u} \|\tilde{y}_{t+1}-y_t\|^2 + \frac{\eta_t\lambda L^2_f}{\rho_l} \|x_{t+1}-x_t\|^2 \nonumber \\
    & \quad + \frac{\eta_t\lambda}{\rho_l}\|\nabla_y f(x_t,y_t)-v_t\|^2 + \frac{L_f\eta^2_t}{2}\|\tilde{y}_{t+1}-y_t\|^2.
\end{align}

Next, using $L_f$-smoothness of function $f(\cdot,y_t)$, such that
\begin{align}
    f(x_t,y_t) + \langle \nabla_x f(x_t,y_t), x_{t+1}-x_t \rangle - \frac{L_f}{2}\|x_{t+1}-x_t\|^2 \leq f(x_{t+1},y_t),
\end{align}
then we have
\begin{align}
    & f(x_t,y_t) -f(x_{t+1},y_t)  \nonumber \\
    & \leq -\langle \nabla_x f(x_t,y_t), x_{t+1}-x_t \rangle + \frac{L_f}{2}\|x_{t+1}-x_t\|^2 \nonumber \\
    & = -\eta_t\langle \nabla_x f(x_t,y_t) - \nabla F(x_t), \tilde{x}_{t+1}-x_t \rangle - \eta_t\langle \nabla F(x_t), \tilde{x}_{t+1}-x_t \rangle + \frac{L_f\eta^2_t}{2}\|\tilde{x}_{t+1}-x_t\|^2 \nonumber \\
    & \leq \frac{\eta_t}{8\gamma}\|\tilde{x}_{t+1}-x_t\|^2 + 2\eta_t\gamma\|\nabla_x f(x_t,y_t) - \nabla F(x_t)\|^2  - \eta_t\langle \nabla F(x_t), \tilde{x}_{t+1}-x_t \rangle + \frac{L_f\eta^2_t}{2}\|\tilde{x}_{t+1}-x_t\|^2 \nonumber \\
    & \leq \frac{\eta_t}{8\gamma}\|\tilde{x}_{t+1}-x_t\|^2 + 2L^2_f\eta_t\gamma \|y_t - y^*(x_t)\|^2 + F(x_t) - F(x_{t+1}) \nonumber \\
    & \quad + \frac{\eta^2_tL}{2}\|\tilde{x}_{t+1}-x_t\|^2 + \frac{\eta^2_tL_f}{2}\|\tilde{x}_{t+1}-x_t\|^2 \nonumber \\
    & \leq \frac{4L^2_f\eta_t\gamma}{\mu} \big(F(x_t) -f(x_t,y_t)\big) + F(x_t) - F(x_{t+1}) + \eta_t(\frac{1}{8\gamma}+\eta_tL)\|\tilde{x}_{t+1}-x_t\|^2,
\end{align}
where the second last inequality is due to Lemma~\ref{lem:A1}, i.e.,
$L$-smoothness of function $F(x)$, and
the the last inequality holds by Lemma~\ref{lem:A2} and $L_f\leq L$.
Then we have
\begin{align} \label{eq:B65}
    F(x_{t+1}) - f(x_{t+1},y_t) & = F(x_{t+1}) - F(x_t) + F(x_t)- f(x_t,y_t) + f(x_t,y_t) -f(x_{t+1},y_t) \nonumber \\
    & \leq (1+\frac{4L^2_f\eta_t\gamma}{\mu}) \big(F(x_t) -f(x_t,y_t)\big) + \eta_t(\frac{1}{8\gamma}+\eta_tL)\|\tilde{x}_{t+1}-x_t\|^2.
\end{align}

Substituting \eqref{eq:B65} into \eqref{eq:B64}, we get
\begin{align}
    & F(x_{t+1}) - f(x_{t+1},y_{t+1})  \\
    & \leq (1-\frac{\eta_t\lambda\mu}{\rho_u})(1+\frac{4L^2_f\eta_t\gamma}{\mu}) \big(F(x_t) -f(x_t,y_t)\big) + \eta_t(\frac{1}{8\gamma}+\eta_tL)(1-\frac{\eta_t\lambda\mu}{\rho_u})\|\tilde{x}_{t+1}-x_t\|^2 \nonumber \\
    & \quad -\frac{\eta_t}{2\lambda\rho_u} \|\tilde{y}_{t+1}-y_t\|^2 + \frac{\eta_t\lambda L^2_f}{\rho_l} \|x_{t+1}-x_t\|^2 + \frac{\eta_t\lambda}{\rho_l}\|\nabla_y f(x_t,y_t)-v_t\|^2 + \frac{L_f\eta^2_t}{2}\|\tilde{y}_{t+1}-y_t\|^2 \nonumber \\
    & = (1-\frac{\eta_t\lambda\mu}{\rho_u})(1+\frac{4L^2_f\eta_t\gamma}{\mu}) \big(F(x_t) -f(x_t,y_t)\big) + \eta_t\big(\frac{1}{8\gamma}+\eta_tL-\frac{\eta_t\lambda\mu}{8\gamma\rho_u}-\frac{\eta^2_tL\lambda\mu}{\rho_u} + \frac{\eta^2_tL^2_f\lambda}{\rho_l}\big)\|\tilde{x}_{t+1}-x_t\|^2 \nonumber \\
    & \quad -\frac{\eta_t}{2}\big(\frac{1}{\lambda\rho_u}-L_f\eta_t\big) \|\tilde{y}_{t+1}-y_t\|^2 + \frac{\eta_t\lambda}{\rho_l}\|\nabla_y f(x_t,y_t)-v_t\|^2 \nonumber \\
    & \leq (1-\frac{\eta_t\lambda\mu}{2\rho_u}) \big(F(x_t) -f(x_t,y_t)\big) + \frac{\eta_t}{8\gamma}\|\tilde{x}_{t+1}-x_t\|^2  -\frac{\eta_t}{4\lambda\rho_u}\|\tilde{y}_{t+1}-y_t\|^2 + \frac{\eta_t\lambda}{\rho_l}\|\nabla_y f(x_t,y_t)-v_t\|^2, \nonumber
\end{align}
where the last inequality holds by $L=L_f(1+\frac{\kappa}{2})$,  $\gamma\leq \min\big(\frac{\lambda\mu}{16\rho_uL},\frac{\rho_l\mu}{16\rho_uL^2_f}\big)$ and $\lambda \leq \frac{1}{2\eta_t L_f\rho_u}$ for all $t\geq 1$, i.e.,
\begin{align}
   & \gamma\leq \frac{\lambda\mu}{16\rho_uL} \Rightarrow \lambda\geq \frac{16\rho_uL\gamma}{\mu}=16\rho_u\gamma(\kappa+\frac{\kappa^2}{2})\geq 8\rho_u\kappa^2\gamma \Rightarrow  \frac{\eta_t\lambda\mu}{2\rho_u} \geq \frac{4L^2_f\eta_t\gamma}{\mu} \nonumber \\
   & \gamma\leq \min\big(\frac{\lambda\mu}{16\rho_uL},\frac{\rho_l\mu}{16\rho_uL^2_f}\big) \Rightarrow \frac{\eta_t\lambda\mu}{8\gamma\rho_u}\geq \eta_tL+ \frac{\eta^2_tL^2_f\lambda}{\rho_l},  \nonumber \\
   &\lambda \leq \frac{1}{2\eta_tL_f\rho_u}  \Rightarrow \frac{1}{2\lambda\rho_u} \geq \eta_t L_f, \ \forall t\geq 1.
\end{align}

\end{proof}

\begin{lemma} \label{lem:D2}
 Suppose that the sequence $\big\{x_t,\tilde{x}_t\big\}_{t=1}^T$ be generated from Algorithm~\ref{alg:2}.
 Let $0<\eta_t \leq 1$ and $0< \gamma \leq \frac{\rho}{2L\eta_t}$,
 then we have
 \begin{align}
  \Psi(x_{t+1}) \leq \Psi(x_t) + \frac{4\gamma L^2_f\eta_t}{\mu\rho}\big(F(x_t)-f(x_t,y_t)\big) + \frac{2\gamma\eta_t}{\rho}\|\nabla_x f(x_t,x_t)-w_t\|^2 -\frac{\rho\gamma\eta_t}{2}\|\mathcal{G}(x_t,w_t,\gamma)\|^2,
 \end{align}
where $\Psi(x)=F(x)+\psi(x)$ and $\mathcal{G}(x_t,w_t,\gamma) = \frac{1}{\gamma}(x_t-\tilde{x}_{t+1})$.
\end{lemma}

\begin{proof}
By the line 4 of Algorithm~\ref{alg:2}, we have
\begin{align} \label{eq:DD1}
\tilde{x}_{t+1} = \mathcal{P}^{\gamma}_{\psi(\cdot)}\big(x_t-\gamma A_t^{-1}w_t\big) = \arg\min_{x\in \mathbb{R}^d}\Big\{ \langle w_t, x\rangle + \frac{1}{2\gamma}(x-x_t)^TA_t(x-x_t) + \psi(x)\Big\}.
\end{align}

By the optimality condition of the subproblem~(\ref{eq:DD1}), we have for any $z\in \mathbb{R}^d$
\begin{align}
 \big\langle w_t + \frac{1}{\gamma}A_t(\tilde{x}_{t+1}-x_t) + \nu_{t+1}, z-\tilde{x}_{t+1}\big\rangle \geq 0,
\end{align}
where $\nu_{t+1}\in \partial \psi(\tilde{x}_{t+1})$.

By using the convexity of $\psi(x)$, and
let $z=x_t$, we can obtain
\begin{align}  \label{eq:DD2}
 \langle w_t, x_t - \tilde{x}_{t+1}\rangle & \geq \frac{1}{\gamma}(\tilde{x}_{t+1}-x_t)^TA_t(\tilde{x}_{t+1}-x_t) + \langle \nu_{t+1}, \tilde{x}_{t+1}-x_t\rangle \nonumber \\
 & \geq \frac{\rho}{\gamma}\|\tilde{x}_{t+1}-x_t\|^2 + \psi(\tilde{x}_{t+1})-\psi(x_t),
\end{align}
where the last inequality holds by Assumption~\ref{ass:4}, i.e., $A_t\succeq \rho I_d$.

Let $\mathcal{G}(x_t,w_t,\gamma) = \frac{1}{\gamma}(x_t-\tilde{x}_{t+1})$. According to Lemma~\ref{lem:A1}, i.e., function $F(x)$ is $L$-smooth, we have
\begin{align} \label{eq:DD3}
  F(x_{t+1}) & \leq F(x_t) + \langle \nabla F(x_t), x_{t+1}-x_t\rangle + \frac{L}{2}\|x_{t+1}-x_t\|^2 \nonumber \\
  & = F(x_t) + \eta_t\langle w_t, \tilde{x}_{t+1}-x_t\rangle + \eta_t\gamma \langle \nabla F(x_t) -w_t, \mathcal{G}(x_t,w_t,\gamma)\rangle+ \frac{\gamma^2\eta_t^2L}{2}\|\mathcal{G}(x_t,w_t,\gamma)\|^2 \nonumber \\
  & \mathop{\leq}^{(i)} F(x_t) - \gamma\rho\eta_t \|\mathcal{G}(x_t,w_t,\gamma)\|^2 - \eta_t\psi(\tilde{x}_{t+1}) + \eta_t\psi(x_t) + \eta_t\gamma \langle \nabla F(x_t) - w_t, \mathcal{G}(x_t,w_t,\gamma)\rangle \nonumber \\
  & \quad + \frac{\gamma ^2\eta_t^2L}{2}\|\mathcal{G}(x_t,w_t,\gamma)\|^2 \nonumber \\
  & = F(x_t) - \gamma\rho\eta_t \|\mathcal{G}(x_t,w_t,\gamma)\|^2 - \eta_t\psi(\tilde{x}_{t+1}) - (1-\eta_t)\psi(x_t) + \psi(x_t) \nonumber \\
  & \quad  + \eta_t\gamma \langle \nabla F(x_t) - w_t, \mathcal{G}(x_t,w_t,\gamma)\rangle+ \frac{\gamma ^2\eta_t^2L}{2}\|\mathcal{G}(x_t,w_t,\gamma)\|^2 \nonumber \\
  & \mathop{\leq}^{(ii)} F(x_t) - \gamma\rho\eta_t \|\mathcal{G}(x_t,w_t,\gamma)\|^2 - \psi(x_{t+1}) + \psi(x_t) \nonumber \\
  & \quad  + \eta_t\gamma \langle \nabla F(x_t) - w_t, \mathcal{G}(x_t,w_t,\gamma)\rangle+ \frac{\gamma ^2\eta_t^2L}{2}\|\mathcal{G}(x_t,w_t,\gamma)\|^2 \nonumber \\
  & \mathop{\leq}^{(iii)} F(x_t) - \frac{\eta_t\gamma\rho}{2}\|\mathcal{G}(x_t,w_t,\gamma)\|^2 - \psi(x_{t+1}) + \psi(x_t) + \frac{\eta_t\gamma}{\rho}\|w_t - \nabla F(x_t)\|^2,
\end{align}
where the inequality (i) holds by the above inequality~(\ref{eq:D2}), and the inequality (ii) is due to $x_{t+1} = x_t + \eta_t(\tilde{x}_{t+1}-x_t)$ and the convexity of function $\psi(x)$, i.e., $\psi(x_{t+1})=\psi((1-\eta_t)x_t + \eta_t\tilde{x}_{t+1})\leq (1-\eta_t)\psi(x_t) + \eta_t\psi(\tilde{x}_{t+1})$,
and the last inequality (iii) holds by $0<\gamma\leq \frac{\rho}{2\eta_t L}$ and the following inequality
\begin{align}
 \langle \nabla F(x_t) -w_t, \mathcal{G}(x_t,w_t,\gamma)\rangle &\leq \|w_t - \nabla F(x_t)\|\|\mathcal{G}(x_t,w_t,\gamma)\| \nonumber \\
 & \leq  \frac{1}{\rho}\|w_t - \nabla F(x_t)\|^2 + \frac{\rho}{4}\|\mathcal{G}(x_t,w_t,\gamma)\|^2,
\end{align}
where the above inequality holds by Young inequality.

Considering the bound of the term $\|\nabla F(x_t)-w_t\|^2$, then we have
\begin{align} \label{eq:DD4}
  \|\nabla F(x_t)-w_t\|^2 & = \|\nabla_x f(x_t,y^*(x_t))- \nabla_x f(x_t,y_t)+\nabla_x f(x_t,y_t)- w_t\|^2\nonumber \\
  & \leq 2\|\nabla_x f(x_t,y^*(x_t))- \nabla_x f(x_t,y_t)\|^2 + 2\|\nabla_x f(x_t,x_t)-w_t\|^2
  \nonumber \\
  & \leq 2L^2_f\|y^*(x_t) - y_t\|^2 + 2\|\nabla_x f(x_t,x_t)-w_t\|^2,
\end{align}
where the first inequality is due to the Cauchy-Schwarz inequality and the second is due to Young's inequality.

Let $\Psi(x) = F(x)+\psi(x)$.
By plugging the above inequalities \eqref{eq:DD4} into \eqref{eq:DD3},
we obtain
\begin{align} \label{eq:DD5}
  \Psi(x_{t+1}) &\leq \Psi(x_t) - \frac{\eta_t\gamma\rho}{2}\|\mathcal{G}(x_t,w_t,\gamma)\|^2 + \frac{\eta_t\gamma}{\rho}\|w_t - \nabla F(x_t)\|^2 \\
  & \leq  \Psi(x_t) - \frac{\eta_t\gamma\rho}{2}\|\mathcal{G}(x_t,w_t,\gamma)\|^2 + \frac{2\eta_t\gamma L^2_f}{\rho}\|y^*(x_t) - y_t\|^2 + \frac{2\eta_t\gamma}{\rho}\|\nabla_x f(x_t,x_t)-w_t\|^2   \nonumber \\
  & \leq \Psi(x_t) - \frac{\eta_t\gamma\rho}{2}\|\mathcal{G}(x_t,w_t,\gamma)\|^2 +\frac{4\eta_t\gamma L^2_f}{\mu\rho}\big(F(x_t)-f(x_t,y_t)\big) + \frac{2\eta_t\gamma}{\rho}\|\nabla_x f(x_t,x_t)-w_t\|^2, \nonumber
\end{align}
where the last inequality holds by the above Lemma~\ref{lem:A2} using in $F(x_t)=f(x_t,y^*(x_t))=\max_y f(x_t,y)$ with $y^*(x_t)\in \arg\max f(x_t,y)$.

\end{proof}

\begin{theorem}  \label{th:A2}
(Restatement of Theorem 2)
 Under the above Assumptions~\ref{ass:1}-\ref{ass:5}, in Algorithm~\ref{alg:2}, let $\eta_t=\frac{k}{(m+t)^{1/3}}$ for all $t\geq 0$, $\alpha_{t+1}=c_1\eta_t^2$, $\beta_{t+1}=c_2\eta_t^2$, $m \geq \max\big(2,k^3, (c_1k)^3, (c_2k)^3\big)$, $k>0$, $ c_1 \geq \frac{2}{3k^3} + \frac{9\rho_u L^2_f}{\rho_l\mu^2}$, $c_2 \geq \frac{2}{3k^3} + \frac{9}{4}$,
 $0< \lambda \leq \min\big(\frac{3}{4\sqrt{2}\mu},\frac{m^{1/3}}{2k L_f\rho_u}\big)$ and $0< \gamma \leq \min\Big(\frac{\lambda\mu}{16\rho_uL},\frac{\rho_l\mu}{16\rho_uL^2_f},\frac{m^{1/3}\rho}{2Lk},\frac{\rho}{8L_f},\frac{\rho^2\lambda\mu^2}{9\rho_uL^2_f}\Big)$, we have
\begin{align}
  \frac{1}{T}\sum_{t=1}^T\mathbb{E}\|\mathcal{G}(x_t,\nabla F(x_t),\gamma)\|  \leq \frac{2\sqrt{3G}m^{1/6}}{T^{1/2}} + \frac{2\sqrt{3G}}{T^{1/3}},
\end{align}
where $\Psi(x)=F(x)+\psi(x)$ and $G = \frac{\Psi(x_1) - \Psi^*}{\gamma k\rho} + \frac{9\rho_u L^2_f}{k\lambda\mu^2\rho^2}\Delta_1  + \frac{2\sigma^2m^{1/3}}{k^2\rho^2} + \frac{2k^2(c_1^2+c_2^2)\sigma^2}{\rho^2}\ln(m+T)$ and
 $\Delta_1 = F(x_1) -f(x_1,y_1)$.
\end{theorem}

\begin{proof}
Since $\eta_t$ is decreasing and $m\geq k^3$, we have $\eta_t \leq \eta_0 = \frac{k}{m^{1/3}} \leq 1$ and $\gamma \leq \frac{\rho}{2L\eta_0} =\frac{m^{1/3}\rho}{2Lk} \leq \frac{1}{2L\eta_t}$ for any $t\geq 0$.
Similarly, $\lambda \leq \frac{1}{2\eta_0 L_f\rho_u} = \frac{m^{1/3}}{2k L_f\rho_u} \leq \frac{1}{2\eta_t L_f\rho_u}$ for all $t\geq 1$
Due to $0 < \eta_t \leq 1$ and $m\geq \max\big( (c_1k)^3, (c_2k)^3 \big)$, we have $\alpha_t = c_1\eta_t^2 \leq c_1\eta_t \leq \frac{c_1k}{m^{1/3}}\leq 1$ and $\beta_t = c_2\eta_t^2 \leq c_2\eta_t \leq \frac{c_2k}{m^{1/3}}\leq 1$.

According to the above Lemma~\ref{lem:B1}, we can obtain
\begin{align}
& \frac{1}{\eta_t}\mathbb{E} \|\nabla_{x} f(x_{t+1},y_{t+1}) - w_{t+1}\|^2 - \frac{1}{\eta_{t-1}}\mathbb{E} \|\nabla_{x} f(x_t,y_t) - w_t\|^2  \\
& \leq \big(\frac{1 - \beta_{t+1}}{\eta_t}  -  \frac{1}{\eta_{t-1}}\big)\mathbb{E} \|\nabla_{x} f(x_t,y_t) -w_t\|^2 + 4L^2_f\eta_t\mathbb{E}\big(\|\tilde{x}_{t+1}-x_t\|^2 + \|\tilde{y}_{t+1}-y_t\|^2\big) +  \frac{2\beta_{t+1}^2\sigma^2}{\eta_t} \nonumber \\
& = \big(\frac{1}{\eta_t} -  \frac{1}{\eta_{t-1}} - c_2\eta_t\big)\mathbb{E} \|\nabla_{x} f(x_t,y_t) -w_t\|^2 + 4L^2_f\eta_t\mathbb{E}\big(\|\tilde{x}_{t+1}-x_t\|^2 + \|\tilde{y}_{t+1}-y_t\|^2\big) +  2c^2_2\eta_t^3\sigma^2, \nonumber
\end{align}
where the second inequality is due to $0<\beta_{t+1}\leq 1$ and $\beta_{t+1}=c_2\eta_t^2$.
Similarly, we have
\begin{align}
& \frac{1}{\eta_t}\mathbb{E} \|\nabla_{y} f(x_{t+1},y_{t+1}) - v_{t+1}\|^2
-  \frac{1}{\eta_{t-1}}\mathbb{E} \|\nabla_{y} f(x_t,y_t) - v_t\|^2  \\
& \leq \big(\frac{1}{\eta_t} - \frac{1}{\eta_{t-1}} - c_1\eta_t\big)\mathbb{E} \|\nabla_{y} f(x_t,y_t) -v_t\|^2 + 4L^2_f\eta_t\mathbb{E}\big(\|\tilde{x}_{t+1}-x_t\|^2 + \|\tilde{y}_{t+1}-y_t\|^2\big) +
 2c^2_1\eta_t^3\sigma^2. \nonumber
\end{align}
By $\eta_t = \frac{k}{(m+t)^{1/3}}$, we have
\begin{align}
\frac{1}{\eta_t} - \frac{1}{\eta_{t-1}} &= \frac{1}{k}\big( (m+t)^{\frac{1}{3}} - (m+t-1)^{\frac{1}{3}}\big) \nonumber \\
& \leq \frac{1}{3k(m+t-1)^{2/3}} = \frac{2^{2/3}}{3k\big(2(m+t-1)\big)^{2/3}} \nonumber \\
& \leq \frac{2^{2/3}}{3k(m+t)^{2/3}} = \frac{2^{2/3}}{3k^3}\frac{k^2}{(m+t)^{2/3}}= \frac{2^{2/3}}{3k^3}\eta_t^2 \leq \frac{2}{3k^3}\eta_t,
\end{align}
where the first inequality holds by the concavity of function $f(x)=x^{1/3}$, \emph{i.e.}, $(x+y)^{1/3}\leq x^{1/3} + \frac{y}{3x^{2/3}}$, and
the last inequality is due to $0<\eta_t\leq 1$.

Let $c_2 \geq \frac{2}{3k^3} + \frac{9}{4}$, we have
\begin{align} \label{eq:GG1}
& \frac{1}{\eta_t}\mathbb{E} \|\nabla_{x} f(x_{t+1},y_{t+1}) - w_{t+1}\|^2 - \frac{1}{\eta_{t-1}}\mathbb{E} \|\nabla_{x} f(x_t,y_t) - w_t\|^2  \\
& \leq -\frac{9\eta_t}{4}\mathbb{E} \|\nabla_{x} f(x_t,y_t) -w_t\|^2 + 4L^2_f\eta_t\mathbb{E}\big(\|\tilde{x}_{t+1}-x_t\|^2 + \|\tilde{y}_{t+1}-y_t\|^2\big) + 2c^2_2\eta_t^3\sigma^2. \nonumber
\end{align}
Let $c_1 \geq \frac{2}{3k^3} + \frac{9\rho_u L^2_f}{\rho_l\mu^2}$, we have
\begin{align} \label{eq:GG2}
& \frac{1}{\eta_t}\mathbb{E} \|\nabla_{y} f(x_{t+1},y_{t+1}) - w_{t+1}\|^2
-  \frac{1}{\eta_{t-1}}\mathbb{E} \|\nabla_{y} f(x_t,y_t) - w_t\|^2  \\
& \leq - \frac{9\rho_u L^2_f\eta_t}{\rho_l\mu^2}\mathbb{E} \|\nabla_{y} f(x_t,y_t) -w_t\|^2 + 4L^2_f\eta_t\mathbb{E}\big(\|\tilde{x}_{t+1}-x_t\|^2 + \|\tilde{y}_{t+1}-y_t\|^2\big) +
2c^2_1\eta_t^3\sigma^2. \nonumber
 \end{align}

According to Lemma \ref{lem:C2}, we have
\begin{align} \label{eq:GG3}
\Psi(x_{t+1})-\Psi(x_t) \leq  + \frac{4\gamma L^2_f\eta_t}{\mu\rho}\big(F(x_t)-f(x_t,y_t)\big) + \frac{2\gamma\eta_t}{\rho}\|\nabla_x f(x_t,x_t)-w_t\|^2 -\frac{\rho\gamma\eta_t}{2}\|\mathcal{G}(x_t,w_t,\gamma)\|^2.
\end{align}
According to Lemma \ref{lem:D2}, we have
\begin{align} \label{eq:GG4}
F(x_{t+1}) -f(x_{t+1},y_{t+1}) - \big(F(x_t) -f(x_t,y_t)\big)  & \leq -\frac{\eta_t\lambda\mu}{2\rho_u} \big(F(x_t) -f(x_t,y_t)\big) + \frac{\eta_t}{8\gamma}\|\tilde{x}_{t+1}-x_t\|^2  \nonumber \\
& \quad -\frac{\eta_t}{4\lambda\rho_u}\|\tilde{y}_{t+1}-y_t\|^2 + \frac{\eta_t\lambda}{\rho_l}\|\nabla_y f(x_t,y_t)-v_t\|^2.
\end{align}

Next, we define a useful Lyapunov function, for any $t\geq 1$
\begin{align}
\Phi_t = \mathbb{E}\big[ \Psi(x_t) + \frac{9\rho_u\gamma L^2_f}{\rho\lambda\mu^2}\big(F(x_t) -f(x_t,y_t) \big) + \frac{\gamma}{\rho}\big(\frac{1}{\eta_{t-1}}\|\nabla_{x} f(x_t,y_t)-w_t\|^2 + \frac{1}{\eta_{t-1}}\|\nabla_{y} f(x_t,y_t)-v_t\|^2 \big) \big].
\end{align}
Let $\mathcal{G}(x_t,w_t,\gamma) = \frac{1}{\gamma}(x_t-\tilde{x}_{t+1})$, then we have
\begin{align}
& \Phi_{t+1}- \Phi_t \nonumber \\
& = \Psi(x_{t+1}) - \Psi(x_t) + \frac{9\rho_u\gamma L^2_f}{\rho\lambda\mu^2}
\Big(F(x_{t+1}) -f(x_{t+1},y_{t+1}) - (F(x_t) -f(x_t,y_t)) \Big)
 + \frac{\gamma}{\rho} \big(\frac{1}{\eta_t}\mathbb{E}\|\nabla_{x} f(x_{t+1},y_{t+1})-w_{t+1}\|^2 \nonumber \\
& \quad  -  \frac{1}{\eta_{t-1}}\mathbb{E}\|\nabla_{x} f(x_t,y_t)-w_t\|^2  +  \frac{1}{\eta_t}\mathbb{E}\|\nabla_{y} f(x_{t+1},y_{t+1})-v_{t+1}\|^2
 -  \frac{1}{\eta_{t-1}}\mathbb{E}\|\nabla_{y} f(x_t,y_t)-v_t\|^2 \big) \nonumber \\
& \leq \frac{4\gamma L^2_f\eta_t}{\mu\rho}\big(F(x_t)-f(x_t,y_t)\big) + \frac{2\gamma\eta_t}{\rho}\|\nabla_x f(x_t,x_t)-w_t\|^2 -\frac{\rho\gamma\eta_t}{2}\|\mathcal{G}(x_t,w_t,\gamma)\|^2\nonumber \\
& \quad  + \frac{9\rho_u\gamma L^2_f}{\rho\lambda\mu^2} \bigg( -\frac{\eta_t\lambda\mu}{2\rho_u} \big(F(x_t) -f(x_t,y_t)\big) + \frac{\gamma\eta_t}{8}\|\mathcal{G}(x_t,w_t,\gamma)\|^2  -\frac{\eta_t}{4\lambda\rho_u}\|\tilde{y}_{t+1}-y_t\|^2 + \frac{\eta_t\lambda}{\rho_l}\|\nabla_y f(x_t,y_t)-v_t\|^2 \bigg)  \nonumber \\
& \quad  -\frac{9\gamma\eta_t}{4\rho}\mathbb{E} \|\nabla_{x} f(x_t,y_t) -w_t\|^2 + \frac{4\gamma L^2_f\eta_t}{\rho}\mathbb{E}\big(\gamma^2\|\mathcal{G}(x_t,w_t,\gamma)\|^2 + \|\tilde{y}_{t+1}-y_t\|^2\big) + \frac{2\gamma c^2_2\eta_t^3\sigma^2}{\rho} \nonumber \\
& \quad -\frac{9\rho_u\gamma L^2_f\eta_t}{\rho_l\rho\mu^2} \mathbb{E}\|\nabla_{y} f(x_t,y_t) -v_t\|^2  + \frac{4\gamma L^2_f\eta_{t}}{\rho }\mathbb{E}\big(\gamma^2\|\mathcal{G}(x_t,w_t,\gamma)\|^2 + \|\tilde{y}_{t+1}-y_t\|^2\big) + \frac{2\gamma c^2_1\eta_t^3\sigma^2}{\rho} \nonumber \\
& \leq -\frac{\gamma L_f^2\eta_t}{2\rho\mu}\big(F(x_t) -f(x_t,y_t) \big) - \frac{\gamma\eta_t}{4\rho}\mathbb{E}\|\nabla_{x} f(x_t,y_t) -w_t\|^2 + (c_1^2+c^2_2)\frac{2\gamma \sigma^2\eta_t^3}{\rho }  \nonumber \\
& \quad -\big(\frac{9\gamma L^2_f}{4\rho\lambda^2\mu^2} - \frac{8\gamma L^2_f}{\rho}\big)\eta_t\mathbb{E}\|\tilde{y}_{t+1} - y_t\|^2 - \big( \frac{\gamma\rho}{2} - \frac{8\gamma^3 L^2_f}{\rho } - \frac{9\rho_u\gamma^2 L^2_f}{8\rho\lambda\mu^2}\big)\eta_t\mathbb{E}\|\mathcal{G}(x_t,w_t,\gamma)\|^2 \nonumber \\
& \leq -\frac{\gamma L_f^2\eta_t}{2\rho\mu}\big(F(x_t) -f(x_t,y_t) \big)  -  \frac{\gamma\eta_t}{4\rho}\mathbb{E}\|\nabla_{x} f(x_t,y_t) -w_t\|^2 - \frac{\rho\gamma\eta_t}{4}\mathbb{E}\|\mathcal{G}(x_t,w_t,\gamma)\|^2 + (c_1^2+c^2_2)\frac{2\gamma \sigma^2\eta_t^3}{\rho },
\end{align}
where the first inequality holds by the above inequalities \eqref{eq:GG1}, \eqref{eq:GG2}, \eqref{eq:GG3} and \eqref{eq:GG4};
the last inequality is due to $0< \lambda \leq \frac{3}{4\sqrt{2}\mu}$ and $0< \gamma \leq \big(\frac{\rho}{8L_f},\frac{\rho^2\lambda\mu^2}{9\rho_uL^2_f}\big) $
for all $t\geq 1$.
Thus, we have
\begin{align} \label{eq:GG5}
 & \frac{L_f^2\eta_t}{2\mu\rho^2}\big(F(x_t) -f(x_t,y_t) \big) + \frac{\eta_t}{4\rho^2}\mathbb{E}\|\nabla_{x} f(x_t,y_t) -w_t\|^2 + \frac{\eta_t}{4}\mathbb{E}\|\mathcal{G}(x_t,w_t,\gamma)\|^2 \nonumber \\
 & \leq \frac{\Phi_t - \Phi_{t+1}}{\rho\gamma} + 2(c_1^2+c^2_2) \frac{\sigma^2\eta_t^3}{\rho^2}.
\end{align}

Taking average over $t=1,2,\cdots,T$ on both sides of \eqref{eq:GG5}, we have
\begin{align}
 & \frac{1}{T} \sum_{t=1}^T \Big( \frac{L_f^2\eta_t}{2\mu\rho^2}\big(F(x_t) -f(x_t,y_t) \big) + \frac{\eta_t}{4\rho^2}\mathbb{E}\|\nabla_{x} f(x_t,y_t) -w_t\|^2 + \frac{\eta_t}{4}\mathbb{E}\|\mathcal{G}(x_t,w_t,\gamma)\|^2  \Big) \nonumber \\
 & \leq  \sum_{t=1}^T \frac{\Phi_t - \Phi_{t+1}}{T\gamma\rho} + \frac{2(c_1^2+c^2_2) \sigma^2}{T\rho^2}\sum_{t=1}^T\eta_t^3. \nonumber
\end{align}
Let $\Delta_1 = F(x_1) -f(x_1,y_1)$, we have
\begin{align} \label{eq:GG6}
 \Phi_1 &= F(x_1) + \frac{9\rho_u\gamma L^2_f}{\rho\lambda\mu^2}\big( F(x_1) -f(x_1,y_1)\big) + \frac{\gamma}{\rho}\big(\frac{1}{\eta_0}\|\nabla_{x}f(x_1,y_1)-w_1\|^2 + \frac{1}{\eta_0}\|\nabla_{y} f(x_1,y_1)-v_1\|^2 \big)  \nonumber \\
 & \leq F(x_1) + \frac{9\rho_u\gamma L^2_f}{\rho\lambda\mu^2}\Delta_1  + \frac{2\gamma\sigma^2}{\rho\eta_0},
\end{align}
where the last inequality holds by Assumption \ref{ass:1}.

Since $\eta_t$ is decreasing, i.e., $\eta_T^{-1} \geq \eta_t^{-1}$ for any $0\leq t\leq T$, we have
 \begin{align}
 & \frac{1}{T} \sum_{t=1}^T \mathbb{E}\Big[  \frac{L_f^2}{2\mu\rho^2}\big(F(x_t) -f(x_t,y_t) \big) + \frac{1}{4\rho^2}\|\nabla_{x} f(x_t,y_t) -w_t\|^2 + \frac{1}{4}\|\mathcal{G}(x_t,w_t,\gamma)\|^2 \Big] \nonumber \\
 & \leq  \sum_{t=1}^T \frac{\Phi_t - \Phi_{t+1}}{\eta_TT\gamma\rho} + \frac{2(c_1^2+c^2_2) \sigma^2}{\eta_TT\rho^2}\sum_{t=1}^T\eta_t^3 \nonumber \\
 & \leq \frac{1}{\eta_TT\gamma\rho}\Big(\Psi(x_1) + \frac{9\rho_u\gamma L^2_f}{\rho\lambda\mu^2}\Delta_1  + \frac{2\gamma\sigma^2}{\rho\eta_0} - \Psi^* \Big) + \frac{2(c_1^2+c^2_2) \sigma^2}{\eta_TT\rho^2}\int^T_1\frac{k^3}{m+t}dt \nonumber \\
 & \leq \frac{\Psi(x_1)-\Psi^*}{\eta_TT\gamma\rho} + \frac{9\rho_u L^2_f}{\rho^2\eta_TT\lambda\mu^2}\Delta_1  + \frac{2\sigma^2}{\eta_TT\eta_0\rho^2} + \frac{2k^3(c_1^2+c^2_2) \sigma^2}{\eta_TT\rho^2}\ln(m+T) \nonumber \\
 & = \bigg( \frac{\Psi(x_1) - \Psi^*}{\gamma k\rho} + \frac{9\rho_u L^2_f}{k\lambda\mu^2\rho^2}\Delta_1  + \frac{2\sigma^2m^{1/3}}{k^2\rho^2} + \frac{2k^2(c_1^2+c_2^2)\sigma^2}{\rho^2}\ln(m+T)\bigg) \frac{(m+T)^{1/3}}{T},
\end{align}
where the second inequality holds by the above inequality \eqref{eq:GG6}. Let $G = \frac{\Psi(x_1) - \Psi^*}{\gamma k\rho} + \frac{9\rho_u L^2_f}{k\lambda\mu^2\rho^2}\Delta_1  + \frac{2\sigma^2m^{1/3}}{k^2\rho^2} + \frac{2k^2(c_1^2+c_2^2)\sigma^2}{\rho^2}\ln(m+T)$,
we have
\begin{align} \label{eq:GG7}
 \frac{1}{T} \sum_{t=1}^T \mathbb{E}\big[ \frac{L_f^2}{2\mu\rho^2}\big(F(x_t) -f(x_t,y_t) \big) + \frac{1}{4\rho^2}\|\nabla_{x} f(x_t,y_t) -w_t\|^2 + \frac{1}{4}\|\mathcal{G}(x_t,w_t,\gamma)\|^2 \big]  \leq \frac{G}{T}(m+T)^{1/3}.
\end{align}

Next, we define a useful gradient mapping
$\mathcal{G}(x_t,\nabla F(x_t),\gamma)=\frac{1}{\gamma}(x_t-x^+_{t+1})$, where $x_{t+1}$
is generated from
\begin{align}
x^+_{t+1} = \arg\min_{x\in \mathbb{R}^d}\Big\{ \langle \nabla F(x_t), x\rangle
+ \frac{1}{2\gamma}(x-x_t)^TA_t(x-x_t) + \psi(x)\Big\}  = \mathcal{P}_{\psi(\cdot)}^\gamma\big(x_t-\gamma A_t^{-1}\nabla F(x_t)\big), \nonumber
\end{align}
where $F(x)=f(x,y^*(x))$ with $y^*(x)\in \arg\max_y f(x,y)$.

According to Lemma~\ref{lem:A1}, we have $\nabla F(x) =\nabla_{x} f(x_t,y^*(x_t))$. Then we
obtain
\begin{align} \label{eq:GG8}
 \|\mathcal{G}(x_t,\nabla F(x_t),\gamma)\| & = \|\mathcal{G}(x_t,\nabla F(x_t),\gamma)-\mathcal{G}(x_t,w_t,\gamma)+\mathcal{G}(x_t,w_t,\gamma)\| \nonumber \\
 & \leq \|\mathcal{G}(x_t,\nabla F(x_t),\gamma)-\mathcal{G}(x_t,w_t,\gamma)\| + \|\mathcal{G}(x_t,w_t,\gamma)\| \nonumber \\
 & = \left\|\frac{1}{\gamma}\big(x_t-\mathcal{P}_{\psi(\cdot)}^\gamma(x_t-\gamma A_t^{-1}\nabla F(x_t))\big)-\frac{1}{\gamma}\big(x_t-\mathcal{P}_{\psi(\cdot)}^\gamma(x_t-\gamma A_t^{-1}w_t)\big)\right\| + \|\mathcal{G}(x_t,w_t,\gamma)\| \nonumber \\
 &\mathop{\leq}^{(i)}  \frac{1}{\rho}\|w_t - \nabla F(x_t)\|+ \|\mathcal{G}(x_t,w_t,\gamma)\| \nonumber \\
 & \leq \frac{1}{\rho}\|w_t - \nabla_{x} f(x_t,y_t)\| +\frac{1}{\rho} \|\nabla_{x} f(x_t,y_t) - \nabla_{x} f(x_t,y^*(x_t)) \|+ \|\mathcal{G}(x_t,w_t,\gamma)\| \nonumber \\
 & \leq \frac{1}{\rho}\|w_t - \nabla_{x} f(x_t,y_t)\| + \frac{L_f}{\rho}\|y_t - y^*(x_t) \|+ \|\mathcal{G}(x_t,w_t,\gamma)\| \nonumber \\
 & \leq\frac{1}{\rho} \|w_t - \nabla_{x} f(x_t,y_t)\| + \frac{\sqrt{2}L_f}{\rho\sqrt{\mu}}\sqrt{F(x_t)-f(x_t,y_t)} + \|\mathcal{G}(x_t,w_t,\gamma)\|,
\end{align}
where the inequality (i) holds by the lemma 2 of~\cite{ghadimi2016mini} and $A_t\succeq \rho I_d$, and the last inequality holds by Lemma~\ref{lem:A2}.

According to the above inequalities~(\ref{eq:GG7}) and~(\ref{eq:GG8}), we have
\begin{align}
 & \frac{1}{T}\sum_{t=1}^T\mathbb{E}\|\mathcal{G}(x_t,\nabla F(x_t),\gamma)\| \nonumber \\
 & \leq \frac{1}{T}\sum_{t=1}^T
\mathbb{E}\big[\frac{1}{\rho}\|w_t - \nabla_{x} f(x_t,y_t)\| + \frac{\sqrt{2}L_f}{\rho\sqrt{\mu}}\sqrt{F(x_t)-f(x_t,y_t)} + \|\mathcal{G}(x_t,w_t,\gamma)\|\big] \nonumber \\
 & \mathop{\leq}^{(i)} \Big( \frac{1}{T} \sum_{t=1}^T \mathbb{E}\big[ \frac{6L_f^2}{\rho^2\mu}\big(F(x_t)-f(x_t,y_t)\big) + \frac{3}{\rho^2}\|\nabla_{x} f(x_t,y_t) -w_t\|^2 +3\|\mathcal{G}(x_t,w_t,\gamma)\|^2 \big] \Big)^{1/2} \nonumber \\
 & \leq \frac{2\sqrt{3G}}{\sqrt{T}}(m+T)^{1/6} \leq \frac{2\sqrt{3G}m^{1/6}}{T^{1/2}} +
 \frac{2\sqrt{3G}}{T^{1/3}},
\end{align}
where the inequality (i) holds by Jensen's inequality.

\end{proof}

\end{onecolumn}

\end{document}